\documentclass{cmslatex}
\usepackage[paperwidth=7in, paperheight=10in, margin=.875in]{geometry}
\usepackage{tikz}
\usepackage[toc,page]{appendix}
\usetikzlibrary{positioning, shapes.geometric}
\usepackage[colorlinks,linkcolor=blue,anchorcolor=green,citecolor=red]{hyperref}

\usepackage{amsfonts,amssymb}
\usepackage{subfigure}
\usepackage{algorithm,algorithmic}
\usepackage{amsmath}
\usepackage{graphicx}
\usepackage{cite}
\usepackage{enumerate}
\renewcommand{\theequation}{\arabic{section}.\arabic{equation}}

   \allowdisplaybreaks
\begin{document}
 \title{{Uncertainty Quantification of Phase Transition Problems with an  Injection Boundary} \thanks{Received date, and accepted date (The correct 
 dates will be enteblack by the editor).}}

            
          \author{Zhenyi Zhang\thanks{LMAM and School of Mathematical Sciences, Peking University, Beijing, 100871, P. R. China. E-mail: \href{mailto:zhenyizhang@stu.pku.edu.cn}{zhenyizhang@stu.pku.edu.cn}.}
          \and Shengbo Ma\thanks{School of Mathematical Sciences, Zhejiang University,
          Hangzhou, 310027, P. R. China. E-mail: \href{mailto:22235054@zju.edu.cn}{22235054@zju.edu.cn}.}
          \and Zhennan Zhou\thanks{Corresponding author.
          Beijing International Center for Mathematical Research, Peking University, Beijing, 100871, 
          P. R. China. E-mail: \href{mailto:zhennan@bicmr.pku.edu.cn}{zhennan@bicmr.pku.edu.cn}.}
          }

         \pagestyle{myheadings} \markboth{
          \uppercase{UQ of phase transition problems with an Injection Boundary}}{Z. ZHANG, S. MA AND Z. ZHOU} \maketitle

          \begin{abstract}
                 {We develop an enthalpy-based modeling and computational framework to quantify uncertainty in Stefan problems with an injection boundary. Inspired by airfoil icing studies, we consider a system featuring an injection boundary inducing domain changes and a free boundary separating phases, resulting in two types of moving boundaries. Our proposed enthalpy-based formulation seamlessly integrates thermal diffusion across the domain with energy fluxes at the boundaries, addressing a modified injection condition for boundary movement. Uncertainty then stems from random variations in the injection boundary. The primary focus of our Uncertainty Quantification (UQ) centers on investigating the effects of uncertainty on free boundary propagation. Through mapping to a reference domain, we derive an enthalpy-based numerical scheme tailored to the transformed coordinate system, facilitating a simple and efficient simulation. Numerical and UQ studies in one and two dimensions validate the proposed model and the extended enthalpy method. They offer intriguing insights into ice accretion and other multiphysics processes involving phase transitions.} 
          \end{abstract}

       \section{Introduction and Motivation}
      Phase transition problems play a pivotal role in understanding fundamental physical processes like melting and solidification. These phenomena have substantial ramifications across various scientific disciplines, including materials science and aerospace engineering (cf. \cite{Vuik93somehistorical,alexiades1992mathematical,book-Stefan-UQ}).
Specifically, these problems involve the intricate dynamics of heat transfer occurring at the interfaces between different phases of materials. This heat transfer, in turn, induces the motion of these phase boundaries. Quantifying and predicting such interplay between heat diffusion and interface movement constitute the fundamental essence of phase transition problems.
And these problems are also commonly referred to as the Stefan problems.

       Phase transition phenomena, coupled with material injection, are also ubiquitous problems in various scientific and engineering problems. One illustrative example motivating research in this field is the challenge posed by ice formation on aircraft airfoils (see, cf. \cite{book-Stefan-UQ,MESSINGER,Myers1,Myers2,Myers3,Myers4shape,Myers5shape,UQflight,Onedime}).
When aircraft operate at high altitudes, they encounter water droplets present in clouds and mist, which can adhere to the aircraft's surfaces, leading to ice formation. Accurately predicting the type and location of ice accumulation is essential for ensuring flight safety and stability. Consequently, this necessity has spurred extensive research into phase transition problems involving injection boundaries, accompanied by corresponding numerical simulations (see, cf. \cite{Myers1, Myers2, Myers3}).

Moreover, quantifying uncertainties in such phase transition problems with injection boundaries has substantial value. Uncertainty stems from influential uncertain factors affecting the dynamics and subsequent dynamic performance. Conducting uncertainty quantification can shed light on the extent to which uncertain conditions impact the problem and foster a deeper understanding of the underlying mechanisms.  However, this remains challenging partially due to the complex interactions between the moving phase interface and injection dynamics in high-dimension. To enable reliable uncertainty analysis, specialized modeling, and computational techniques are required to account for these coupled physics.

In this work, we aim to develop an integrated framework combining UQ techniques with phase transition modeling and simulations tailored for the phase transition problems with an injection boundary. 
The enthalpy-based model 
we considered in two-dimension in this paper is
\begin{align}
              &\rho c \frac{\partial U}{\partial t}=k \Delta U, \quad (y,x,t)\in (0, 1)\times (0, S(y,t)) \times(0, \infty), \\
      &\rho c \frac{\partial U}{\partial t}=k \Delta U,  \quad (y,x,t)\in (0, 1)\times (S(y,t),L(t)) \times(0, \infty), \\
       &U(y,S(y,t), t)=T_{m}, \label{eq:Stefan1} \\
      &\rho L_h{{v_n}}(y,t)(t)=k \left[\nabla U\cdot 
 \boldsymbol{n}\right](y,S(y,t), t), \label{eq:Stefan2}  \\
      &  k \frac{\partial U}{\partial x}(y,L(t),t)+\frac{d L(t)}{dt} H(y,L(t),t)=\eta(y,t)  \frac{d L(t)}{dt}.
      \label{eq:IBC}
      \end{align}
      The spatial domain representing the region occupied by the materials in both phases is  \((y,x)\in [0,1]\times [0, L(t)]\), where $y$ is the horizontal position and $x$ is the vertical position. In this context,  $U$ represents the temperature and $H$ represents the enthalpy (or the total heat energy, to be introduced in Sec. \ref{sec:model}) with $\rho, c$ denoting the density and specific heat capacity respectively.  The governing system describes heat diffusion across two phases and incorporates two key moving boundaries: the internal free boundary and the external injection boundary.

The internal dynamic interface $S(y,t)$ separates the solid and liquid phases. At the interface, the temperature equals a critical value, denoted by $T_m$, as shown in Eq. \eqref{eq:Stefan1}. The moving normal velocity ${v_n}$ of the free boundary is governed by the latent heat $L_h$ during phase transition, along with the heat fluxes differences across the interface, as in Eq. \eqref{eq:Stefan2}. Here, $\boldsymbol{n}$ represents the unit normal vector at the interface and $[g]$ denotes the difference of $g$ with respect to the two-side limits of the interface. We refer to Eq. \eqref{eq:Stefan1} and Eq. \eqref{eq:Stefan2} as the Stefan condition. 
 
\begin{figure}[htbp]
        \centering
        \includegraphics[width=8cm]{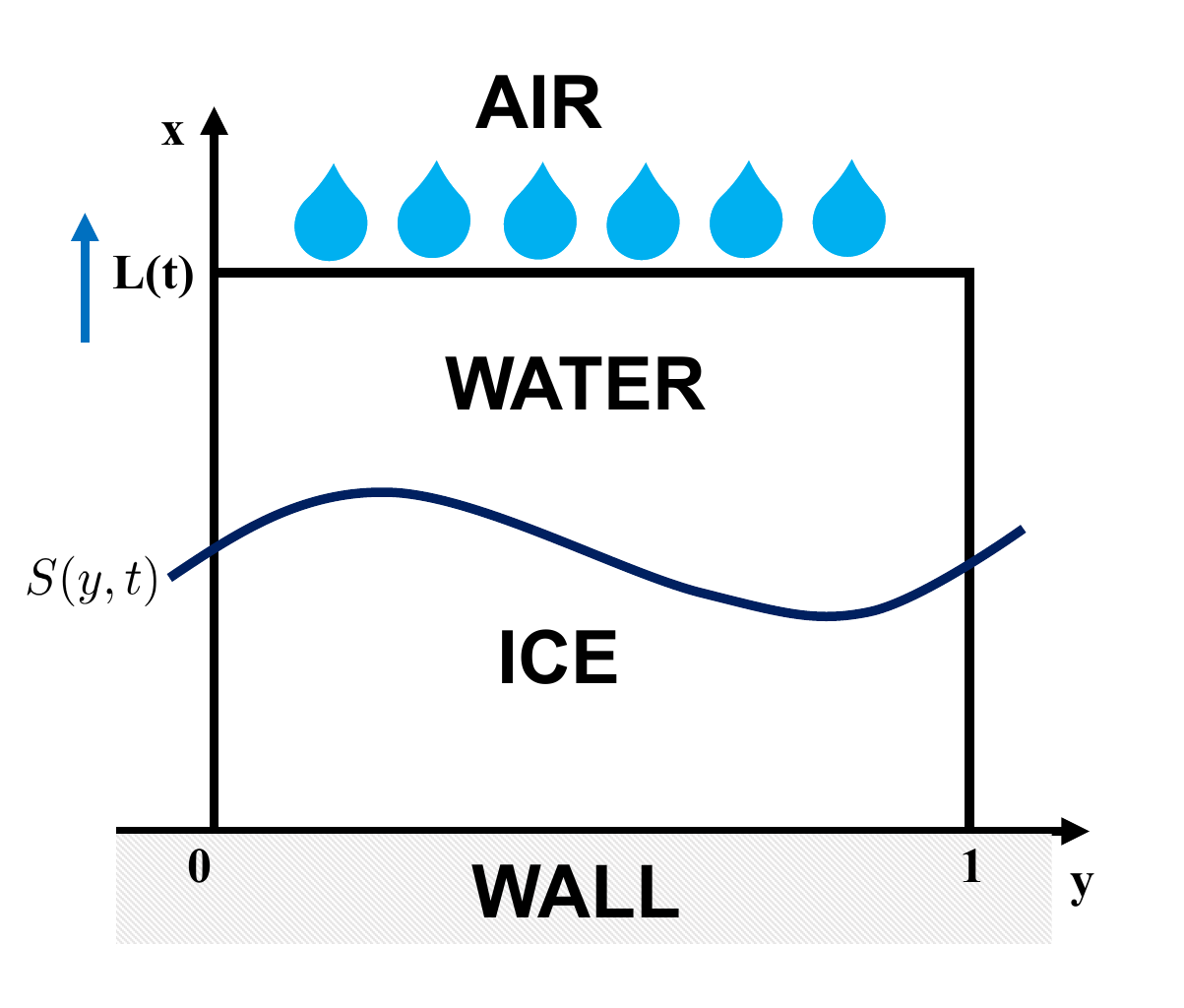}
        \caption{Model system for the 2D phase transition problem with an injection boundary.}
        \label{introfig:Twodsche}
      \end{figure}

The external injection boundary, denoting the inlet of incoming mass, is represented by $x=L(t)$. The heat conductivity is $k$ and $\eta(y,t)$ denotes the spatial-dependent heat parameter accounting for the energy influx carried by incoming mass flux. A key modeling design is to treat $L(t)$ as a homogeneously growing boundary to reflect the even distribution of injected water (see, Fig. \ref{introfig:Twodsche}). Yet we allow flexibility for $\eta(y,t)$ to capture potential spatial energy inhomogeneity. Accordingly, appropriate energy-conserving boundary conditions are imposed on $L(t)$, as described in Eq. \eqref{eq:IBC}.
However, ensuring appropriate treatment of the boundary condition at the injection boundary poses nontrivial challenges. As one of the main results of this work, we put forward a thermodynamically consistent adjustment to tackle such difficulties intrinsically connected to the injection boundary dynamics, which will be carefully designed and analyzed below.

Our proposed enthalpy-based model shares some commonalities but also differs in some key aspects from the previous models in icing formation.
 After Stefan formulated the first mathematical formulation for the liquid ice two-phase problem, Messenger (cf. \cite{MESSINGER}) and Meyers (see, cf. \cite{Myers1,Myers2,Myers3,Myers4shape,Myers5shape}) successively developed a formulation which explains the two different mechanisms associated with rime and glaze ice formation. These icing models, while useful and computationally efficient, have neglected important factors in the physical process.
First, the Myers and Messenger models are 1D models, but they have been applied to simulate the 3D aircraft icing (see, cf. \cite{Myers4shape,Myers5shape,Onedime}). Therefore, the integrated model cannot capture horizontal thermal diffusion and the appropriate coupling of convection and diffusion in the icing process.
Second, the modification of the boundary condition caused by the movement at the injection boundary has not been considered.
Most notably, the boundary condition at the injection boundary in the previous models only consisted of a heat flux term (see, cf. \cite{MESSINGER,Myers1,Myers2,Myers3,Myers4shape,Myers5shape}), without incorporating terms related to injection boundary movement velocity, i.e., $\frac{d L(t)}{dt} H(y,L(t),t)$ as in Eq. (\ref{eq:IBC}).

       The study of uncertainty quantification for the proposed phase transition model is of great importance.
        However, the movement of both the injection and phase boundaries poses challenges in designing suitable and robust uncertainty quantification algorithms.
 To address this, in this work, we aim to develop a computational framework integrating UQ for phase transition with injection boundaries. The key components include:
      \begin{itemize}
          \item An enthalpy-based model incorporating thermal diffusion and incoming droplet energy transfer. This provides an accurate mathematical representation.
          \item An extended enthalpy numerical method implicitly handling the coupled physics and moving boundaries. This enables efficient solutions tailored to injection problems.
          \item  Customized non-intrusive UQ techniques adapted to the enthalpy formulation and injection physics, including specialized sampling, gPC, and sensitivity analysis. This allows quantifying uncertainties in parameters like injection rate.
      \end{itemize}
      The proposed integrated framework with targeted UQ advances computational modeling capabilities for quantifying uncertainties in multi-phase problems involving interacting injection and phase transition.  
      
      The primary challenge in designing algorithms within this framework lies in effectively solving the proposed phase transition model.
      While numerical algorithms have been constructed and studied for similar problems, most are restricted to 1D scenarios.
Front tracking methods with fixed grids require complex front tracking algorithms and can struggle with irregular interface geometries in higher dimensions (cf. \cite{Myers2,VSG1}). Coordinate transform techniques developed for one-phase problems do not directly extend to two-phase systems (cf. \cite{TwoMoving}).
An alternative is the enthalpy method, which implicitly handles phase transition interfaces by reformulating the governing equations to automatically satisfy free boundary conditions. This approach has been successfully leveraged to model phase transition in a wide range of contexts including the Stefan problem (see, cf. \cite{StefanEnthalpy1,StefanEnthalpy2-Date,StefanEnthalpy3-Voller,StefanEnthalpy4,StefanEnthalpy5}, 
etc.)  and convection-diffusion scenarios (cf. 
\cite{Conv-diffEnthalpy1,Conv-diffEnthalpy2,Conv-diffEnthalpy3}). However, the application of enthalpy methods to phase transition problems with injection boundaries has been lacking thus far.

The key focus of this paper is to bridge the limitation mentioned earlier by incorporating some key contributions:
\begin{itemize}
          \item We develop a modified injection boundary condition that conserves the thermal energy of impacting droplets by seamlessly combining the inlet heat flux with boundary growth.
          \item We incorporate multidimensional thermal diffusion, significantly enhancing model fidelity and capture accuracy compared to vertical-only diffusion.
          \item  We extend the enthalpy formulation to handle two-dimensional physics with an injection boundary. This yields an efficient simulation framework without tracking explicit interfaces.
      \end{itemize}   
Collectively, these meaningful improvements are numerically validated. Custom UQ techniques then enable uncertainty analysis.

            The rest of the paper is organized as follows.
In Section \ref{sec:model}, we derive a model incorporating the coupled physics in detail.
 Section \ref{sec:comput} develops the moving boundary enthalpy solver.
 In Section \ref{sec:UQ}, we adapt non-intrusive UQ techniques to the injection physics.
                  Section \ref{sec:Numer} provides 1D and 2D numerical studies validating our integrated framework.
Finally, we summarize our findings and discuss
promising future directions in Section \ref{sec:Conclu}.

\emph{Notations and symbols} ---
We consider a spatial domain \(\Omega \subset \mathbb{R}^d\),
where \(d=1,2\) in our discussion.
And \(\Omega\) can be divided into three disjoint regions. We denote 
by \(A\) the liquid region, 
\(B\) the solid region.
 \(S\) is the common portion of the 
boundaries of \(A\) 
and \(B\). And \(D=A  \cup B \) ,  \(\Omega=A\cup B \cup S\).
We also define the space-time cylinder set 
\(Q_T=\Omega \times [0,T]\). 
$Q_T$ can also be partitioned into three regions.
The spatial and temporal normal vector to $S\times [0,T]$ at $(\boldsymbol{x}, t)$ is  $\boldsymbol{\nu}=\left(\boldsymbol{\nu}_{\boldsymbol{x}}, v_n\right)$.
The space \(C_{0}^{\infty}(Q_T)\) contains 
smooth functions with compact
support in \(Q_T\).
The subscripts \(s\) and \(l\) represent the solid and liquid phases, respectively.

          \section{Model Derivation based on Enthalpy}\label{sec:model}

Phase transition problems with injection boundaries pose modeling challenges due to moving boundary interactions. Previous models make simplifying assumptions about the injection boundary condition.
Our objective is to derive a model coupling thermal diffusion and influx dynamics across the domain and through the injection boundary. 

 The core concept in our model derivation is the so-called enthalpy. Enthalpy represents the total energy, incorporating both sensible heat and latent heat, where the former causes the temperature change and the latter governs the phase transition. 
         
\begin{figure}[htbp]
    \centering
    \includegraphics[width=9cm]{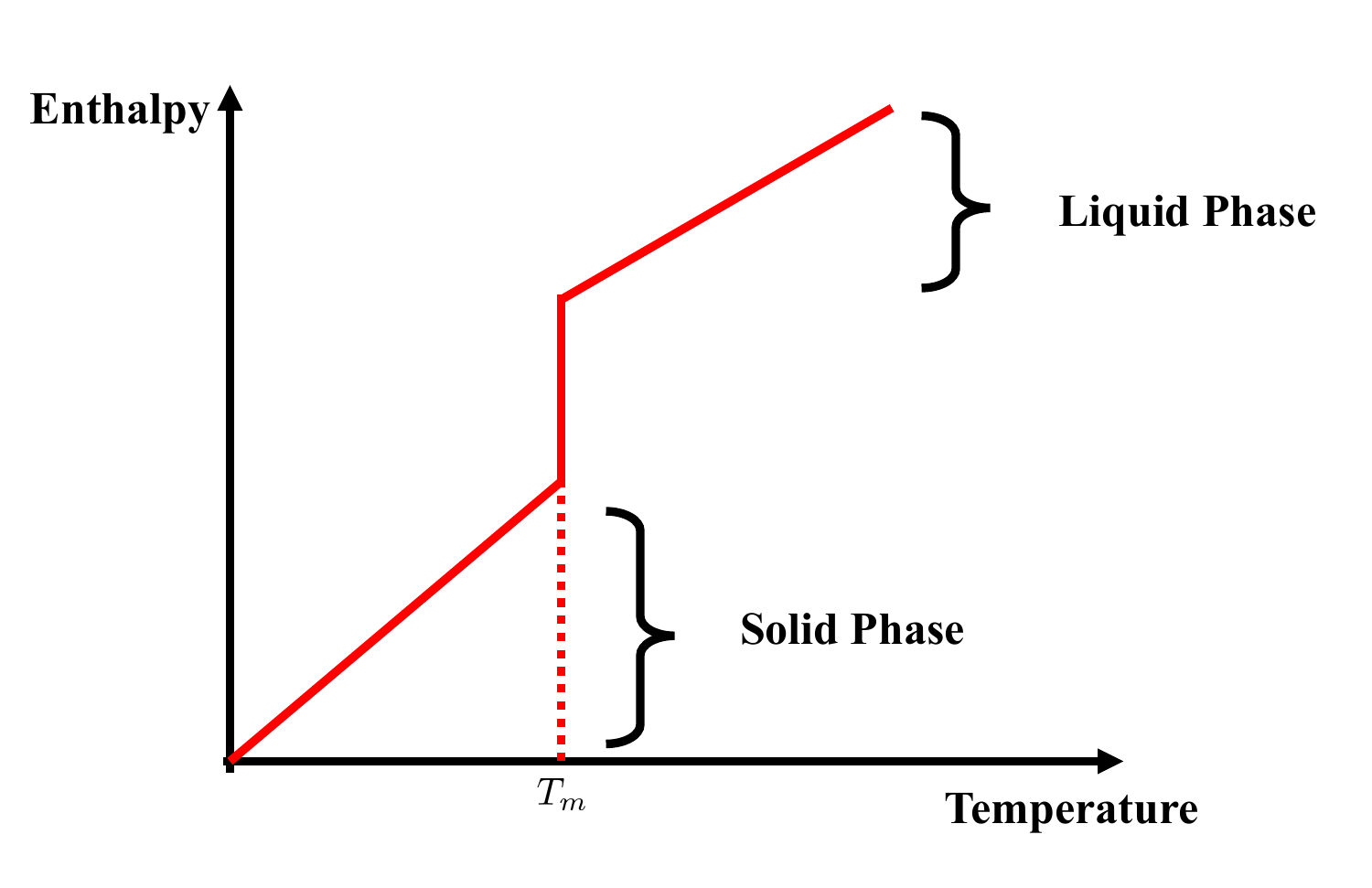}
    \caption{Enthalpy curve.  The curve depicts the enthalpy change as a function of temperature. Below the critical temperature $T_m$, the curve reflects the enthalpy change of the solid phase growths gradually with temperature. At $T_m$, a sharp discontinuity marks the absorption of latent heat during the solid-to-liquid phase transition. Beyond $T_m$, the enthalpy of the liquid phase exhibits a similar ascent with respect to temperature increments.}
    \label{fig:Enthalpy}
  \end{figure}

          More specifically, if we denote enthalpy by \(H(t)\) and it is
 defined as follows:
   \begin{equation}
     H(t)=
   \begin{cases}
       c\rho U(t), \quad U<T_m\\
       c\rho U(t)+\rho L_h, \quad U>T_m.
   \end{cases}
   \label{Enthalpy}
   \end{equation}
  Here, we recall that  \(U(t)\), \(\rho\), \(c\)  and $L_h$ are the temperature,
   density, specific heat capacity and latent heat, respectively.
   \(T_m\) represents the critical phase transition temperature
   and we set  \(T_m=0\) 
   in our subsequent discussion.
   In Fig. \ref{fig:Enthalpy}, we depict the enthalpy curve. Note that due to the discontinuity of enthalpy, piecewise integration is also required when performing spatial integration. 
  
We introduce the basic model assumptions as follows.

     \medskip
 
   \begin{assumption}\label{assu:energy} In the model system, there are two phases of the material.
   \begin{itemize} 
       \item In the interior of each phase, the evolution of the heat distribution is governed by the heat equation.
       \item   The dynamics of the phase transition curve are governed by the Stefan condition Eq. \eqref{eq:Stefan1} and Eq. \eqref{eq:Stefan2}.
       \item  At the injection boundary, the energy influx is proportional to the incoming mass flux, where the proportional coefficient $\eta(y,t)$ can be spatially and temporally dependent. 
   \end{itemize}
   \end{assumption}

     \medskip
     
 Assumption \ref{assu:energy} is made to simplify the exposition. We can also account for heat transfer with air and consider more detailed energy conservation, such as energy dissipation, evaporation, etc. (see e.g. \cite{Myers2,Myers3}). 
Also for simplicity, we assume the density $\rho$, specific heat capacity $c$, and thermal conductivity $k$ have the same values in the solid and liquid phases. The overall results we develop here extend as well to account for different options,  at the cost of more technical calculations with no essential difficulty. 

   We summarize our model as follows
\begin{equation}
\begin{cases}
\begin{aligned}
    &  \left.
    \begin{aligned}
        & \rho c \frac{\partial U}{\partial t}=k \Delta U,  \quad \text { in }(0, 1)\times (0, S(y,t)) \times(0, \infty), \\
        & \rho c \frac{\partial U}{\partial t}=k \Delta U, \quad \text { in }(0, 1)\times (S(y,t),L(t)) \times(0, \infty),
    \end{aligned}\right\}
   (\text{Governing PDEs}) \\
    & \left.
    \begin{aligned}
      & U(y,S(y,t), t)=T_{m}, \\
        & \rho L_h{{v_n}}(y,t)(t)=k \left[\nabla U\cdot 
 \boldsymbol{n}\right](y,S(y,t), t),\\
        
    \end{aligned}
    \right\}(\text{Stefan Condition})\\   
     & \left.
     \begin{aligned}
    & k \frac{\partial U}{\partial x}(y,L(t),t)+\frac{d L(t)}{dt} H(y,L(t),t)=\eta(y,t)  \frac{d L(t)}{dt}.
    \end{aligned}
    \right\}(\text{Injection B. C.})\\
\end{aligned}
\end{cases}
\label{eq:category}
\end{equation}
   The system is completed with the initial condition $U(\cdot,\cdot, 0)=T_{\text{\text{initial}}}$ and the fixed boundary conditions:  $U(y, 0, \cdot)=T_0, \text { on }(0, \infty)$
      and periodic boundary conditions in $y$.  

In view of Assumption \ref{assu:energy}, the governing equation and the Stefan condition are natural components of the model. Therefore, the model derivation boils down to justifying the injection boundary condition. Adhering to the energy conservation law, we perform an enthalpy balance across the global spatial domain to derive the injection boundary condition. The changes in the total energy stem from both heat transfer from outer boundaries, while the evolution of the interface of the phases preserves energy. 

         \subsection{Deriving Injection B. C. in Two-dimensional Models}
  Recall that the spatial domain is the 
     rectangular  region \((y,x)\in [0,1]\times [0, L(t)]\),
     where \(y\) is the horizontal spatial variable and \(x\)
     is the vertical spatial variable. Note that the form of the spatial domain implies an additional model assumption.

     \medskip
     
     \begin{assumption}\label{assu:TwoDIncrease}
  We assume the injection boundary grows homogeneous in $y$.
     \end{assumption}

     \medskip
    
     As we mentioned earlier, the Assumption \ref{assu:TwoDIncrease} indicates that the growth of the injection boundary is in the form of a flat layer (see, Fig. \ref{introfig:Twodsche}). 
      This assumption is reasonable at least for a short period in the applications of interests.
      Our computational framework may be extended to irregular boundary growth cases, but this requires substantial mathematical treatment and goes beyond the scope of this paper. However, recall that in our modeling, the energy brought by injection is inhomogeneous in spatial and time.

 Next, we embark on the derivation process, showcasing the step-by-step procedure to establish the injection boundary condition based on energy conservation.

 For the two-dimension model, $$
 \frac{d}{dt} \iint_{[0,1]\times [0, L(t)]}H(y,x,t)dydx$$ indicates the change of the total enthalpy within this region, which is contributed by outer boundary conditions.
 
For one thing,
\begin{equation}
    \frac{d}{dt} \iint_{[0,1]\times [0,L(t)]}H(y,x,t)dydx= \frac{d}{dt} \left[\iint_{A}H(y,x,t)dydx+\iint_{B}H(y,x,t)dydx\right],
\end{equation}
where $A$ and $B$ represent the liquid region and the solid region respectively. Note that $A$ has two moving boundaries and  $B$ has one moving boundary.

Thus by applying the Leibniz integral rule, we obtain
\begin{equation}
 \frac{d}{dt} \iint_{A}H(y,x,t)dydx= \iint_{A}H_t(y,x,t)dydx+\int_{S}H(y,x,t)v_nds
 +\int_{0}^{1}H(y,L(t),t)\frac{d L(t)}{dt}dy.
\end{equation}
With $H_t=c\rho U_t$ and the heat equation $c\rho U_t=k\Delta U$ both in $A$ and $B$,  we obtain
\begin{equation}\label{eq2.6}
    \begin{aligned}
    \frac{d}{dt} \iint_{A}H(y,x,t)dydx&= \iint_{A}k \Delta Udydx+\int_{S}H(y,x,t)v_nds
 +\int_{0}^{1}H(y,L(t),t)\frac{d L(t)}{dt}dy \\
 &= \int_{\partial A}k \frac{\partial U}{\partial \boldsymbol{n}}ds+\int_{S}H(y,x,t)v_nds+\int_{0}^{1}H(y,L(t),t)\frac{d L(t)}{dt}dy\\
 &=\int_{S} k \nabla U \cdot \boldsymbol{n }ds+\int_{\partial A \cap \partial Q_T} k \nabla U \cdot \boldsymbol{n }ds+\int_{S}H(y,x,t)v_nds\\
 &+\int_{0}^{1}H(y,L(t),t)\frac{d L(t)}{dt}dy,\\
    \end{aligned}
\end{equation}
\noindent where in the second equation, we apply the integration by parts formula.
Note that in $A$ in addition to the movement of the free boundary, there is also a movement of the injection boundary, while only the movement of the free boundary in $B$.

Similarly in $B$ we have, 
\begin{equation}\label{eq2.7}
    \begin{aligned}
     \frac{d}{dt} \iint_{B}H(y,x,t)dydx&= \iint_{B}k \Delta Udydx+\int_{S}H(y,x,t)v_nds\\
 &= \int_{\partial B}k \frac{\partial U}{\partial \boldsymbol{n}}ds+\int_{S}H(y,x,t)v_nds\\
 &=-\int_{S} k \nabla U \cdot \boldsymbol{n }ds+\int_{\partial B \cap \partial Q_T} k \nabla U \cdot \boldsymbol{n }ds-\int_{S}H(y,x,t)v_nds.\\
    \end{aligned}
\end{equation}
\noindent By the definition of $H$, we add Eq. (\ref{eq2.6}) and Eq. (\ref{eq2.7}) , we have
\begin{equation}
    \begin{aligned}
    \frac{d}{dt} \iint_{[0,1]\times [0,L(t)]}H(y,x,t)dydx&= \frac{d}{dt} \left[\iint_{A}H(y,x,t)dydx+\iint_{B}H(y,x,t)dydx\right]\\
    &=\int_{S} \rho L_h v_n+k\left(\nabla U \cdot \boldsymbol{n }-\nabla U \cdot \boldsymbol{n }\right) ds+\int_{\partial Q_T}k \nabla U \cdot \boldsymbol{n }ds\\
    &+\int_{0}^{1}H(y,L(t),t)\frac{d L(t)}{dt}dy.
    \end{aligned}
\end{equation}    
\noindent Combing with Stefan condition:
\begin{equation}
    \begin{aligned}
\rho L_h{{v_n}}(y,t)(t)=k \left[\nabla U\cdot 
 \boldsymbol{n}\right](y,S(y,t), t),
     \end{aligned}
\end{equation}  
Then we can get
\begin{equation}
    \begin{aligned}
    \frac{d}{dt} \iint_{[0,1]\times [0,L(t)]}H(y,x,t)dydx&= \frac{d}{dt} \left[\iint_{A}H(y,x,t)dydx+\iint_{B}H(y,x,t)dydx\right]\\
    &=\int_{\partial Q_T}k \nabla U \cdot \boldsymbol{n }ds+\int_{0}^{1}H(y,L(t),t)\frac{d L(t)}{dt}dy\\
    &=\int_{0}^{1}\left[k\frac{\partial U}{\partial x}(y,L(t),t)+H(y,L(t),t)\frac{d L(t)}{dt} \right] dy\\
    -k&\int_{0}^{1} \frac{\partial U}{\partial x}(y,0,t)dy
    +k\int_{0}^{L(t)}\left[\frac{\partial U}{\partial y}(1,x,t)-\frac{\partial U}{\partial y}(0,x,t)\right] dx.
    \end{aligned}
\end{equation}
Thus, the change in system energy is contributed by airfoil heat transfer, the incoming water droplets, and thermal diffusion at transverse boundaries.
The energy of the water droplets entering can be expressed as
\begin{equation}
 E_l=\int_{0}^{1} \eta(y,t)\frac{d L(t)}{dt} dy.
\end{equation}
It implies that 
\begin{equation}
    k \frac{\partial U}{\partial x}(y,L(t),t)+\frac{d L(t)}{dt} H(y,L(t),t)=\eta(y,t) \frac{d L(t)}{dt}.
\end{equation}
With Assumption \ref{assu:energy} and
Assumption \ref{assu:TwoDIncrease} we can then propose our 
two-dimensional model as it shown in Eq. (\ref{eq:category}).

The boundary condition at the injection boundary (Eq. (\ref{eq:category})) couples the conductive heat flux with the enthalpy flux due to the boundary movement. This differs from previous models that considered only a prescribed heat flux condition without accounting for the dynamics of the moving boundary. The key novelty is incorporating the enthalpy flux term $\frac{d L(t)}{dt} H(y,L(t),t)$ that arises from applying energy conservation at the moving boundary. This provides a more complete characterization of the energy transfer at the interface between the phases.

         \subsection{One-dimensional Model}
         
        As a direct corollary to our two-dimensional model, let \(L(t)\) be the whole length of the domain and $\left[0, L(t)\right]$ is our required region, then we propose our one-dimensional model as follows:

          \begin{figure}[htbp]
            \centering
            \includegraphics[width=6cm]{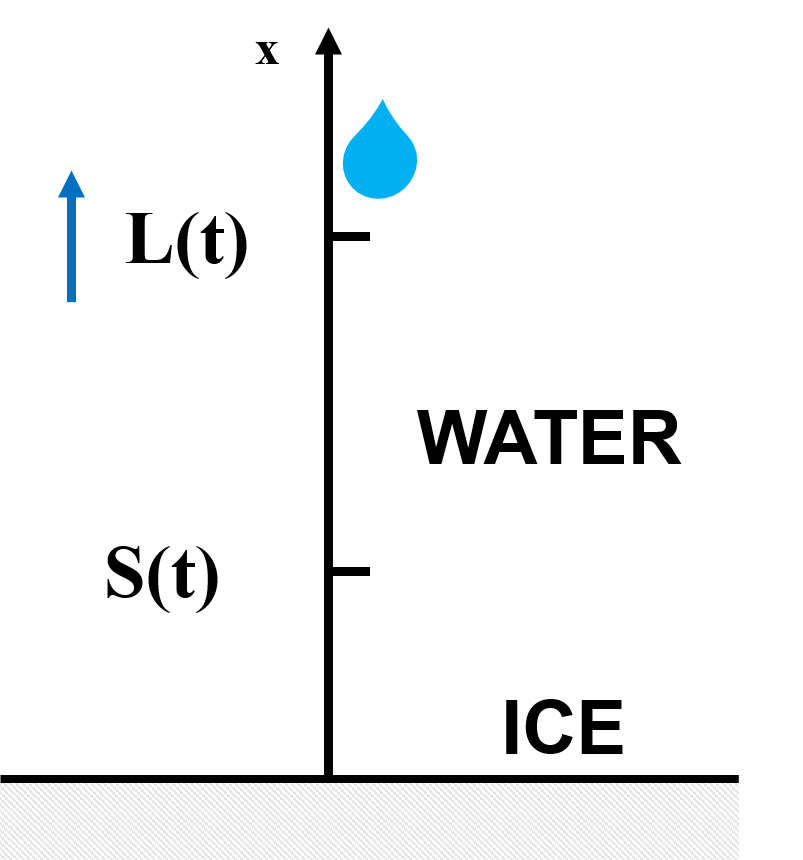}
            \caption{Model system for the 1D phase transition problem with an injection boundary. }
            \label{fig:OneDsch}
          \end{figure}

  \begin{equation}
    \begin{cases}
      \begin{aligned}
          &\rho c \frac{\partial U}{\partial t}=k \frac{\partial^2 U}{\partial x^2}, \quad 0<x<S(t), \quad t>0, \\
  &\rho c \frac{\partial U}{\partial t}=k \frac{\partial^2 U}{\partial x^2}, \quad S(t)<x<L(t), \quad t>0 ,\\
  &U(0, t)={T_0}, \quad t>0, \\
  &U(S(t), t)=T_m, \quad t>0,\\
  &U(x, 0)=T_{\text{initial}},  \quad S(0)=0,\\
  &  k \frac{\partial U}{\partial x}(L(t),t)+\frac{d L(t)}{dt} H(L(t),t)=\eta(t)  \frac{d L(t)}{dt} ,\\
  &\rho L_h \frac{d S}{d t}=\left[k\frac{\partial U}{\partial x}\right], x=S(t), t>0.
      \end{aligned}
  \end{cases}
  \label{model:1}
  \end{equation}
At the initial time, the temperature of the material is 
  uniform, \(T_{\text{initial}}\), higher than the phase-transition 
  temperature 0. 
  The material is cooled at the point \(x = 0 \)
  by imposing a constant temperature \(T_0\), lower than the phase-transition
temperature 0.

   \subsection{Uncertainty Inputs}
   To account for uncertainties, we model some parameters as random fields. Specifically, we take the injection coefficient \(\beta(t,w)\) and influx energy \(\eta(y,t,w)\) (or \(\eta(t,w)\) in one-dimension) to be random processes indexed by  $w$ in probability space. This introduces stochasticity into the boundary dynamics and energy flux. 
   
   A main goal of our uncertainty quantification analysis will be to characterize how this randomness propagates through the coupled multiphysics system to impact quantities of interest like the phase boundary location. The randomness poses challenges for robust numerical simulation and uncertainty propagation. Our specialized sampling and computational methods will be designed to enable reliable quantification of these effects (see, Sec. \ref{sec:comput} and Sec. \ref{sec:UQ}).

      \section{Computational Framework}\label{sec:comput}
     The proposed models in both Eq. (\ref{model:1}) and Eq. (\ref{eq:category}) involve a two-phase Stefan problem with coupled free and injection boundary movements, posing computational challenges. Existing numerical methods exhibit limitations in capturing complex boundary dynamics in high dimensions. To address this, we develop an extended enthalpy method implicitly handling the boundaries. This avoids explicitly tracking interfaces, enabling simpler and more robust treatment in multiple dimensions. The key innovation involves introducing a coordinate transformation mapping the time-dependent domain to a reference domain. For the transformed problem, we derive a new generalized enthalpy governing equation. This extends traditional enthalpy formulations for Stefan problems to incorporate boundary injection.  We first detail the formulation for the one-dimensional model.

         \subsection{One-dimensional Computational Framework}\label{subsec:onecomput}
We develop an extended numerical enthalpy method to solve the proposed model in Eq. (\ref{model:1}). The key steps involve:
\begin{itemize}
    \item Introducing a coordinate transformation and deriving a new governing enthalpy equation incorporating the transformed Stefan condition.
    \item Discretizing the enthalpy equation using implicit finite differences and relating enthalpy and temperature to implicitly track phase boundaries.
\end{itemize}

We will explain each component in detail.
  To transform the moving domain to a fixed reference domain, we first introduce the following coordinate transformation:
\begin{equation}
  z=\frac{x}{L(t)},\quad z\in [0,1], 
\end{equation}
where $L(t)$ is the length of the domain.  And we  naturally define $\ S^{\star}(t)=\frac{S(t)}{L(t)}$.

Under this scaling, and the governing equations in Eq. (\ref{eq:category}) transform to:
\begin{equation}\label{eq:oneConvert}
  \begin{cases}
    \begin{aligned}
        &\rho c\frac{\partial U}{\partial t}=\frac{k}{L^2(t)}\frac{\partial^2 U}{\partial z^2}+\rho c\frac{\partial U}{\partial z}\frac{z}{L(t)}L'(t), \ x \in A, \, B, \ 
 t>0, \\
&U(S^{\star}(t), t)=U(S^{\star}(t), t)=T_m, \quad t>0,\\
& \frac{k}{L(t)} \frac{\partial U}{\partial z}(1,t)+
\frac{d L(t)}{dt} H(1,t)=\eta(t)  \frac{d L(t)}{dt}, \quad t>0 ,\\
&U(z, 0)=T_{\text{initial}}, \quad 0 \leq z \leq 1,\\
&\rho L_h\frac{\partial S^{\star}}{\partial t}=\frac{1}{L(t)}\left[(\frac{k}{L(t)}\frac{\partial U}{\partial z})\right]-\rho L_h\frac{S^{\star}(t)}{L(t)}L'(t),\quad z=S^{\star}(t).  
    \end{aligned}
\end{cases}
\end{equation}

The details of the variable transformations are provided in Appendix \ref{Appendix B}.
With the change of variables, the problem now involves a fixed domain with a single moving interface. This facilitates enthalpy treatment as outlined next.

\begin{figure}[htbp]
  \centering
   \includegraphics[width=10cm]{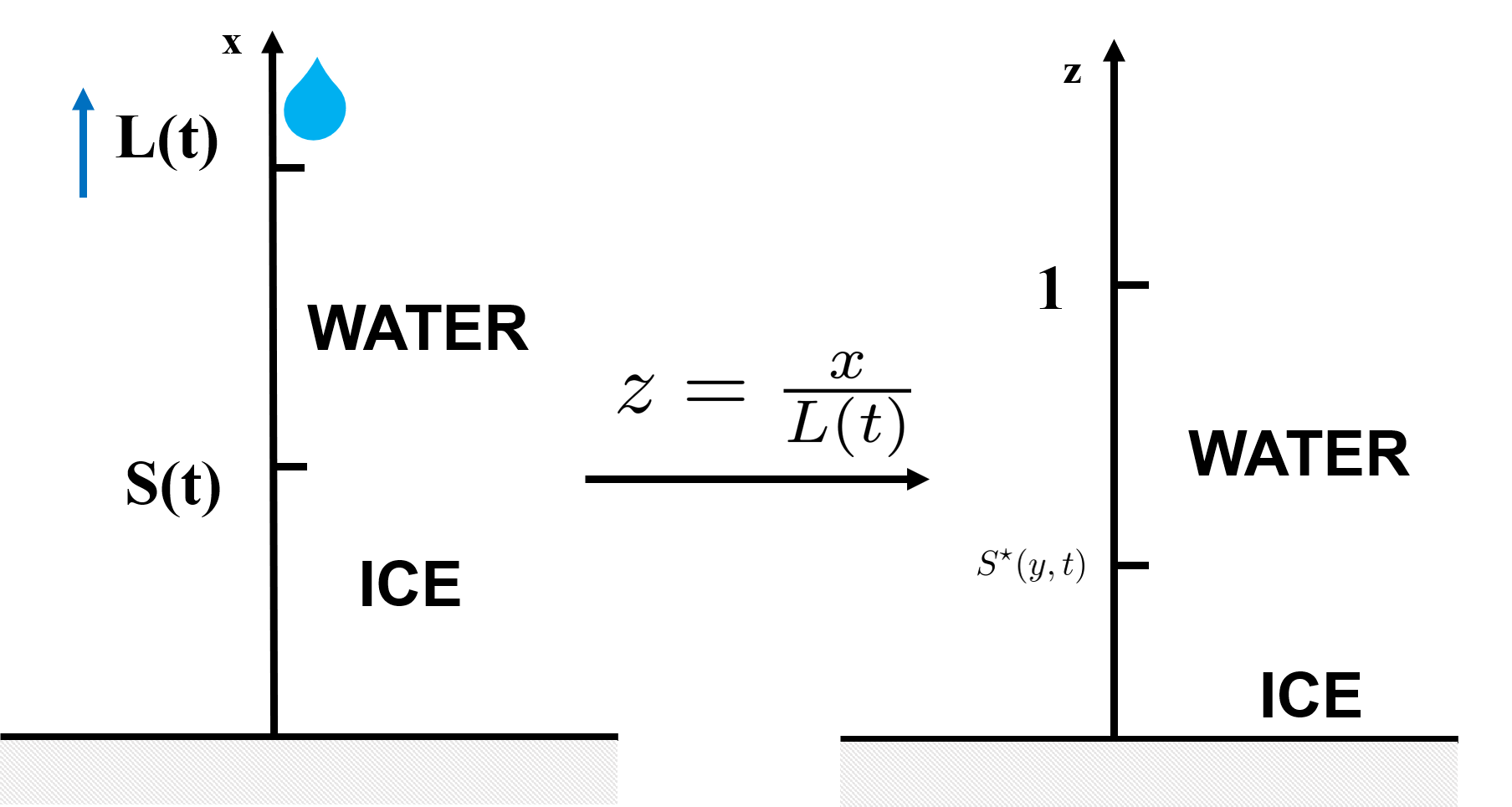}
   \caption{ The transformed system compared to the original elementary system
        for the phase transition problem with an injection boundary in one-dimension.
        The  reference system is
         depicted.  } 
   \label{fig:OneDDtrans} 
 \end{figure}

 The coordinate scaling transforms the original model into a new form with additional convection terms as it is shown in Eq. (\ref{eq:oneConvert}). Unlike traditional Stefan problems, this convection arises from geometric stretching rather than fluid transport. The altered Stefan condition also differs, as scaling introduces an extra convection component. While the transformed system Eq. (\ref{eq:oneConvert}) resembles a convection-diffusion Stefan model (cf. \cite{Conv-diffModel1,Conv-diffModel2,Conv-diffEnthalpy1,Conv-diffEnthalpy2}), key distinctions exist due to the modified convection physics and Stefan condition. 

Next, we introduce the enthalpy method (cf. \cite{StefanEnthalpy2-Date,StefanEnthalpy3-Voller,StefanEnthalpy4}), which reformulates the phase change problem in terms of enthalpy $H$ instead of temperature. Enthalpy represents the total energy and implicitly incorporates latent heat effects.
The key steps of the enthalpy approach are:
\begin{itemize}
    \item Derive a governing equation for the enthalpy field and this equation will have a discontinuity at the phase interface.
\item  Discretize the enthalpy equation using finite differences and solve the discrete system to obtain the enthalpy field.
\item  The interface location corresponds to the zero level set of the enthalpy and temperature can be recovered from the enthalpy values.
\end{itemize}

To solve the transformed problem, we tailor the enthalpy formulation by deriving a governing equation incorporating the altered Stefan condition. As in traditional enthalpy approaches, the key step involves obtaining an equation for the enthalpy field that implicitly handles the moving boundary. For the scaled model, we derive the following governing enthalpy equation for Eq. (\ref{eq:oneConvert}):
\begin{equation}\label{eq:gover1}
  \frac{\partial H}{\partial t}=\frac{k}{L^2(t)}\frac{\partial^2 U}{\partial z^2}+\frac{\partial H}{\partial z}\frac{z}{L(t)}L'(t),
\end{equation}
where
\(H\) is defined in Eq. (\ref{Enthalpy}).
Note that if \(L'(t)=0\), it degenerates to 
the enthalpy equation of traditional Stefan problem.
The detailed derivation of the enthalpy method will be provided in Appendix \ref{Appendix A.}.

Now we just need to solve the enthalpy governing equation
Eq. (\ref{eq:gover1}).
Next we design the numerical scheme for solving Eq. (\ref{eq:gover1}).
Since  $z$ has been dimensionless, we next  define a reference length $L_{\text{ref}}$, and use it to nondimensionalize the rest of  variables.
To facilitate understanding, we  introduce the following dimensionless  variables:
\begin{equation}
  \begin{aligned}
    \varphi=\frac{H}{\rho L_h},~ 
    \theta=\frac{CU}{L_h},~ 
    Z=z,~ 
    \tau=\frac{k}{\rho C L_{\text{ref}}^2}t, ~
    \end{aligned}
\end{equation}
where $\varphi$, $\theta$, $\tau$ are the nondimensionalized enthalpy, temperature, and time, respectively. And we take $\theta_{0}$, $\theta_{\text{initial}}$, $\theta_m$ to represent the nondimensionalization of $T_{0}$, $T_{\text{initial}}$, $T_m$.

We assume that \(L(\tau)=1+\widehat{\beta}\tau\), where \(\widehat{\beta}=\frac{\rho C L_{\text{ref}}^2}{k}\beta\). 
Let $\gamma=\frac{k L_h}{C L_{\text{ref}}^2}$ and
$\widetilde{\eta}=\frac{k}{\rho C L_{\text{ref}}^2}\eta $, We set $\widehat{\eta} = \frac{\widetilde{\eta}}{\gamma}$ as a dimensionless parameter which reflects the incoming heat influx (IHI) carried by the injection mass at the boundary,  then the boundary condition at \(L(t)\) can be written as 
\begin{equation}
  \frac{L_{\text{ref}}^2}{L(\tau)}\frac{\partial \theta}{\partial z}(1,\tau)+ \widehat{\beta}  \varphi(1,\tau) =\widehat{\eta}   \widehat{\beta}. 
\end{equation}

Then the Eq. (\ref{eq:oneConvert})
can be converted into the following form:
\begin{equation}\label{eq:oneUndimeConvert}
  \begin{cases}
    \begin{aligned}
        &  \frac{\partial \theta}{\partial \tau}=
        \frac{L_{\text{ref}}^2}{L^2(\tau)}
        \frac{\partial^2 \theta}{\partial z^2}+
        \frac{\partial \theta}{\partial z}\widehat{\beta}\frac{z}{ L(\tau)},  \ x \in A, \, B, \ 
 \tau>0,  \\
&\theta(0, \tau)=\theta_{0}, \quad \tau>0, \\
&\theta(S^{\star}(\tau), \tau)=\theta(S^{\star}(\tau), \tau)=\theta_m, \quad \tau>0,\\
&   \frac{L_{\text{ref}}^2 }{L(\tau)}\frac{\partial \theta}{\partial z}(1,\tau)+\widehat{\beta}  \varphi(1,\tau) =\widehat{\eta}   \widehat{\beta}, 
\quad \tau>0 ,\\
&\theta(z, 0)=\theta_{\text{initial}}, \quad 0 \leq z \leq 1,\\
&\frac{\partial S^{\star}}{\partial \tau}=
\frac{L_{\text{ref}}^2 }{L(\tau)}\left[\frac{1}{L(\tau)}\frac{\partial \theta}{\partial z}\right]-
\frac{S^{\star}(\tau)}{L(\tau)}\widehat{\beta}, z=S^{\star}(\tau).
    \end{aligned}
\end{cases}
\end{equation}

From our discussion above, it is known that
 only the enthalpy equation needs to be solved, and our enthalpy equation Eq. (\ref{eq:gover1}) 
 now yields the following form:

\begin{equation}\label{eq:govercon1}
  \frac{\partial \varphi}{\partial \tau}=\frac{L_{\text{ref}}^2}{L^2(\tau)}\frac{\partial^2 \theta}{\partial z^2}+\frac{\partial \varphi}{\partial z}\widehat{\beta}\frac{z}{ L(\tau)},
\end{equation}
with 
\begin{equation}
  \begin{aligned}
    \theta&=\varphi, \text{ for } \varphi<0 \text{   (solid)  }, \\
    \theta&=0,  \text{ for } 0<\varphi<1 \text{   (phase transition)  }, \\
    \theta&=\varphi-1, \text{ for } \varphi\geq 1 \text{   (liquid)  }.
\end{aligned}
\end{equation}

The above  relations  can be generalized as
\begin{equation}
\label{relat}
  \theta=\varphi+\widetilde{\varphi},\quad
\mbox{where} \quad 
  \widetilde{\varphi}=0.5(|1-\varphi|-|\varphi|-1).
\end{equation}
It is important to note that $\widetilde{\varphi}=0$ in solid, $\widetilde{\varphi}=-\varphi$ during phase transition, and $\widetilde{\varphi}=-1.0$ in liquid.
This will bring more conveniences in tracking the location
of the inner interface.

Next, we construct the numerical scheme. Let $[0,T]$ be the time evolution region
and $\Delta \tau$ be the time step size.
We introduce a uniform time grid $\tau^{(1)}=0\le \tau^{(2)}\le \cdots \le \tau^{(\frac{T}{\Delta \tau}+1)}=T$, where $\tau^{(k)}=(k-1)\Delta \tau$, $k=1,\, ...\, , \, \frac{T}{\Delta \tau}+1:=K$.
Similarly we choose the computational domain as $[0,1]$ and a uniform mesh grid $z_1=0\le z_2 \le \cdots \le z_{(\frac{1}{\Delta z}+1)}=1$ with uniform grid spacing $\Delta z$,  where $z_i=(i-1)\Delta z, \, i=1,\,\cdots \,,\,\frac{1}{\Delta z}+1:=N$.
We apply the implicit finite-difference scheme
for the diffusion term and the upwind scheme for the convection term
to discretize Eq. (\ref{eq:govercon1}), i.e., 
\begin{equation}
  \frac{\varphi_i^{(k+1)}-\varphi_i^{(k)}}{\Delta \tau}=\frac{L_{\text{ref}}^2}{L^2(\tau)}\frac{\theta_{i-1}^{(k+1)}-2 
  \theta_i^{(k+1)}+\theta_{i+1}^{(k+1)}}{(\Delta z)^2}+\widehat{\beta} 
  \frac{1}{ L(\tau)}z_{i} \frac{\varphi_{i+1}^{(k+1)}-\varphi_{i}^{(k+1)}}{\Delta z}.
\end{equation}

The discretization of the boundary at \(z=1\) is 
\begin{equation}
  \frac{L_{\text{ref}}^2}{L(\tau)}\frac{\theta_{N}^{(k+1)}-\theta_{N-1}^{(k+1)}}{\Delta 
  z}+\widehat{\beta}  \varphi_N^{(k+1)}=\widehat{\eta} \widehat{\beta}.
\end{equation}

Utilize the relation (\ref{relat}), we obtain that 
\begin{equation}
\begin{aligned}
    &\left(1+2\frac{D L_{\text{ref}}^2}{L^2(\tau)}+\frac{\widehat{\beta} z_i}{L(\tau)}
    \Lambda\right) \varphi_i^{(k+1)}
    -\frac{D L_{\text{ref}}^2}{L^2(\tau)}
    \varphi_{i-1}^{(k+1)}-\left(\frac{D L_{\text{ref}}^2}{L^2(\tau)}+\frac{\widehat{\beta}z_i}{L(\tau)}
    \Lambda \right)\varphi_{i+1}^{(k+1)}\\
    &=\frac{D L_{\text{ref}}^2}{L^2(\tau)}\left(\widetilde{\varphi}_{i-1}^{(k+1)}-2  \widetilde{\varphi}_i^{(k+1)}+\widetilde{\varphi}_{i+1}^{(k+1)}\right)+\varphi_{i}^{(k)}, \,  i= 2, \cdots, N-1,
\end{aligned}
\end{equation}
in the interior of the domain, where $\Lambda= \frac{\Delta \tau}{\Delta z}$ and $D= \frac{\Delta \tau}{(\Delta z)^2}$.
And at the boundary \(z=1\), it follows that:
\begin{equation}
    \left(\frac{L_{\text{ref}}^2}{L(\tau) \Delta z}+\widehat{\beta}\right)\varphi_N^{(k+1)}
    -\frac{L_{\text{ref}}^2}{L(\tau) \Delta z} \varphi_{N-1}^{(k+1)} =\frac{L_{\text{ref}}^2}{L(\tau) \Delta z}\left(\widetilde{\varphi}_{N-1}^{(k+1)}- 
\widetilde{\varphi}_N^{(k+1)}\right)+
   \widehat{\eta}\widehat{\beta}.
\end{equation}

In the discretized equations, the right hand side contains \(\widetilde{\varphi}^{(k+1)}\) terms representing the phase change interface position updated at the current iteration. As these depend on the latest $\varphi$ solution, we approximate \(\widetilde{\varphi}^{(k+1)}\) using the previous $\varphi^{(k)}$values. 
This approximation allows discretizing the enthalpy equation without having to solve additional coupled equations for the interface position. The resulting equations contain only $\varphi$ terms, enabling a simple tridiagonal matrix solution
 and it can be solved by applying
 the tridiagonal matrix algorithm (TDMA).

\subsection{ Two-dimensional Computational Framework}

\begin{figure}[htbp]
  \centering
  \includegraphics[width=13cm]{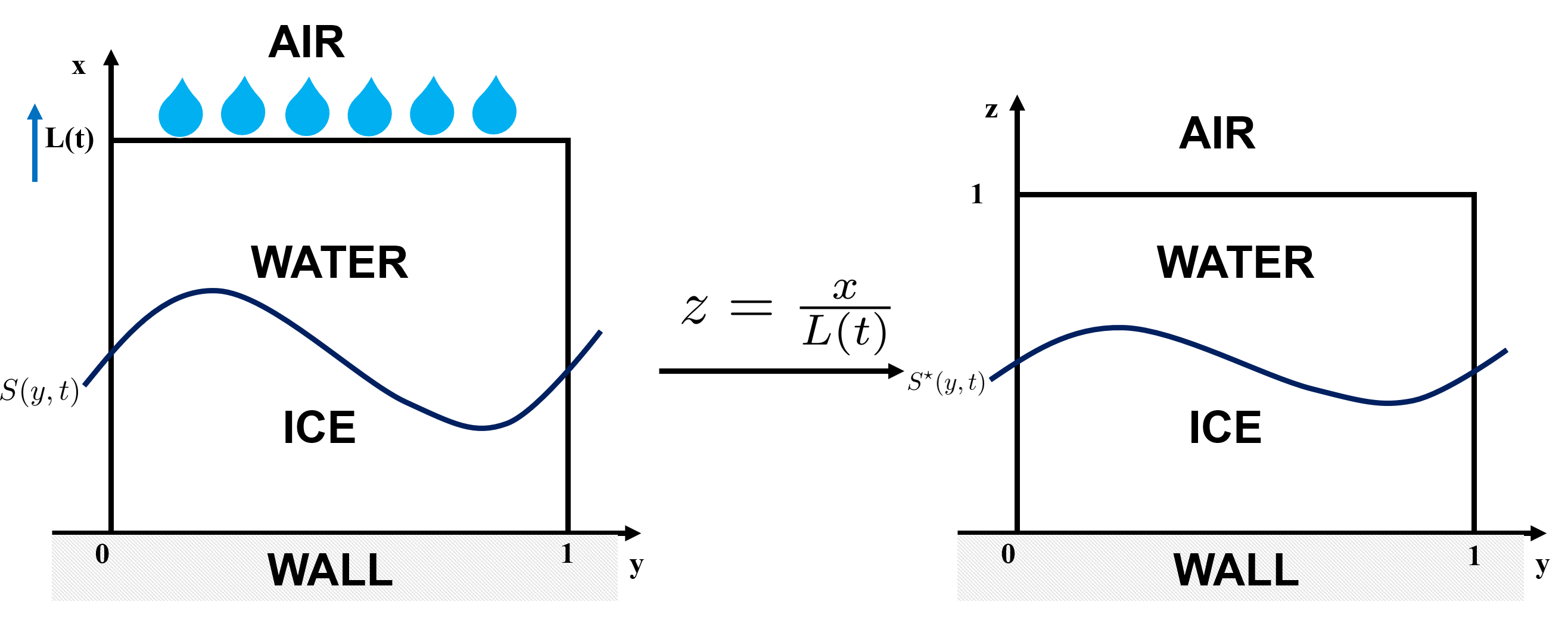}
  \caption{  The transformed system compared to the original elementary system
        for the phase transition problem with an injection boundary in two-dimension.
        The  reference system is
         depicted.   } 
  \label{fig:TWODDtrans} 
\end{figure}

Similar to the one-dimensional case, in our two-dimensional model in 
Eq. (\ref{eq:category}), we also apply the coordinate transformation, but 
in this case we only  make changes in the vertical \(x\) direction.
We take 
\begin{equation}
    z=\frac{x}{L(t)},
\end{equation}
and 
we write $u(y,x,t)=U(y,z(x,t),t)$, where \(u(y,x,t)\) represents
\(U(y,x,t)\) in the Eq. (\ref{eq:category}), as shown in Fig. \ref{fig:TWODDtrans}.

Then the Eq. (\ref{eq:category}) can be converted as follows:
\begin{equation}\label{eq:TwoDmodel}
  \begin{cases}
    \begin{aligned}
        &\rho c\frac{\partial U}{\partial t}=\frac{k}{L^2(t)}\frac{\partial^2 U}{\partial z^2}+k\frac{\partial^2 U}{\partial y^2}+\rho c\frac{\partial U}{\partial z}\frac{z}{L(t)}L'(t), \, 0<z<\frac{S(y,t)}{L(t)}:=S^{\star}(y,t), \\
&\rho c\frac{\partial U}{\partial t}=\frac{k}{L^2(t)}\frac{\partial^2 U}{\partial z^2}+k\frac{\partial^2 U}{\partial y^2}+\rho c\frac{\partial U}{\partial z}\frac{z}{L(t)}L'(t), \, S^{\star}(y,t)<z<1, \\
&U(y, 0, \cdot)=T_0, \quad \text { on }(0, \infty) , \\
&U(y,S^{\star}(y,t), t)=T_{m},\\
&\frac{k}{L(t)} \frac{\partial U}{\partial z}(1,t)+\frac{d L(t)}{d t} H(1, t)=\eta \frac{d L(t)}{d t}, \quad t>0 ,\\
&U(\cdot,\cdot, 0)=T_{\text{initial}}, \\
&\rho L_h{{v_z}}(y,t)(t)=\frac{1}{L(t)}\left[\frac{k}{L(t)}\frac{\partial U}{\partial z}\right]-\rho L\frac{S^{\star}(y, t)}{L(t)}L'(t), \, z=S^{\star}(y, t).
    \end{aligned}
\end{cases}
\end{equation}
Similarly, we apply the enthalpy method to solve, 
and the governing equation is
\begin{equation}
  \frac{\partial H}{\partial t}=\frac{k}{L^2(t)}\frac{\partial^2 U}{\partial z^2}+k\frac{\partial^2 U}{\partial y^2}+\frac{\partial H}{\partial z}\frac{z}{L(t)}L'(t).
\end{equation}
The discussion of the governing equation is almost same as in the one-dimensional case,
hence we omit it here.
By introducing the same variable as in one-dimensional case, 
the Eq. (\ref{eq:TwoDmodel}) now reads that:
\begin{equation}\label{eq:TwoUndimeConvert}
  \begin{cases}
    \begin{aligned}
        &  \frac{\partial \theta}{\partial \tau}=
        \frac{L_{\text{ref}}^2}{L^2(\tau)}
        \frac{\partial^2 \theta}{\partial z^2}+
        L_{\text{ref}}^2\frac{\partial^2 \theta}{\partial y^2}+
        \frac{\partial \theta}{\partial z}\widehat{\beta}\frac{z}{ L(\tau)}, \, 0<z<\frac{S(y,\tau)}{L(\tau)}:=S^{\star}(y,\tau), \\
&\frac{\partial \theta}{\partial \tau}=
\frac{L_{\text{ref}}^2}{L^2(\tau)}
\frac{\partial^2 \theta}{\partial z^2}+ L_{\text{ref}}^2 \frac{\partial^2 \theta}{\partial y^2}+
\frac{\partial \theta}{\partial z}\widehat{\beta}\frac{z}{ L(\tau)}, \, S^{\star}(y,\tau)<z<1, \\
&\theta(y, 0, \cdot)=\theta_0 \quad \text { on }(0, \infty) ,\\
&\theta(y,S^{\star}(y,\tau), \tau)=\theta_{m},\\
&   \frac{ L_{\text{ref}}^2 }{L(\tau)}\frac{\partial \theta}{\partial z}(1,\tau)+\widehat{\beta} \varphi(1,\tau) =\widehat{\eta}  \widehat{\beta}, 
\quad \tau>0 ,\\
&\theta(\cdot,\cdot, 0)=\theta_{\text{initial}},
\\
&{{v_z}}(y,\tau)(\tau)=
\frac{L_{\text{ref}}^2}{L(\tau)}\left[\frac{1}{L(\tau)}
\frac{\partial \theta}{\partial z}\right]-
\frac{S^{\star}(y, \tau)}{L(\tau)}\widehat{\beta}, \, z=S^{\star}(y, \tau).
    \end{aligned}
\end{cases}
\end{equation}

The enthalpy governing equation
now admits the following form:

\begin{equation}
  \frac{\partial \varphi}{\partial \tau}=\frac{L_{\text{ref}}^2}{L^2(\tau)}
  \frac{\partial^2 \theta}{\partial z^2}+
  L_{\text{ref}}^2\frac{\partial^2 \theta}{\partial y^2}+
  \frac{\partial \varphi}{\partial z}\widehat{\beta} \frac{z}{ L(\tau)}.
\end{equation}
In the two-dimensional case, we can also apply the implicit finite difference
method to discretize, for example the alternating direct implicit method or
the predictor-corrector method, etc.
But for simplicity, here we use the explict finite difference method to 
discretize since it is easy to implement.
With a slight abuse of notation, we first introduce the space-time grid.
Let $[0,T]$ be the time evolution region
and $\Delta \tau$ be the time step size.
We introduce a uniform time grid $\tau^{(1)}=0\le \tau^{(2)}\le \cdots \le \tau^{(K)}=T$, where $\tau^{(k)}=(k-1)\Delta \tau, k=1,...,K=\frac{T}{\Delta \tau}+1$.
Similarly we choose the computational domain as $[0,1]\times[0,1]$ and a uniform mesh grid with  uniform grid spacing $\Delta y$ and $\Delta y$ respectively,  $\mathcal{T}=\{(y_j,z_i), y_j=(j-1)\Delta y, z_i=(i-1)\Delta z, j=1,\cdots,N_y=\frac{1}{\Delta y}+1, i=1,\cdots,N=\frac{1}{\Delta z}+1\}$.
Then the scheme takes the form:
\begin{equation}
  \begin{aligned}
    \frac{\varphi_{j,i}^{(k+1)}-
    \varphi_{j,i}^{(k)}}{\Delta \tau}&=
    \frac{L_{\text{ref}}^2}{L^2(\tau)}
    \frac{\theta_{j,i-1}^{(k)}-2 \theta_{j,i}^{(k)}+
    \theta_{j,i+1}^{(k)}}{(\Delta z)^2}+
    \frac{\widehat{\beta}}{L(\tau)}z_{i} 
    \frac{\varphi_{j,i+1}^{(k)}-\varphi_{j,i}^{(k)}}{\Delta z}\\
    &+L_{\text{ref}}^2\frac{\theta_{j-1,i}^{(k)}-2 \theta_{j,i}^{(k)}+\theta_{j+1,i}^{(k)}}{(\Delta y)^2}.
  \end{aligned}
\end{equation}
And at the horizontal boundary \(z=1\), we have:
\begin{equation}
  \theta_{j,N+1}^{(k)}=(1-\frac{L(\tau)}{ L_{\text{ref}}^2 }\Delta z \widehat{\beta})\theta_{j,N}^{(k)}+
  \widehat{\beta}(\widehat{\eta}-1)\frac{L(\tau)}{ L_{\text{ref}}^2 }\Delta z, j=1, \dots,N_y.
\end{equation}

Under some assumptions, the convergence of the enthalpy method can be 
proved (cf.\cite{StefanEnthalpy1,StefanEnthalpy2-Date,StefanEnthalpy3-Voller,StefanEnthalpy4,StefanEnthalpy5}). With the solver, we can now study 
the uncertainty quantification of parameters.

\section{Uncertainty Quantification}\label{sec:UQ}
In this section, we will develop specialized techniques to quantify uncertainties in the phase transition model with injection boundary derived in Sec. \ref{sec:model}.
Recall our definitions of $\widehat{\eta}$ and $\widehat{\beta}$ in Sec. \ref{sec:comput}. The non-dimensionalized parameter $\widehat{\eta}$ represents the incoming heat influx (IHI) carried by the injection mass at the boundary. The IHI $\widehat{\eta}$ quantifies the enthalpy that must be attained at the boundary for the incoming water to be in equilibrium and it must be greater than 1 since the water incoming.
The parameter $\widehat{\beta}$ represents the injection boundary velocity (IBV) governing the rate of mass influx.

We focus on the impact of uncertainty in the IBV  $\widehat{\beta}$ and the corresponding IHI  $\widehat{\eta}$.
Let $\xi$ denote the random parameters (e.g., $\widehat{\beta}$, $\widehat{\eta}$). We assume $\xi$ follows a known probability distribution $p(\xi)$ and
$S(t,\xi)$ represents the free boundary location, which now depends on the random parameters $p(\xi)$.
Then $\{\xi_m\}^M_{m=1}$ denotes a set of $M$ samples generated from the distribution of $\xi$ and
$\phi_n(\xi)$ are orthogonal polynomial basis functions used in the generalized polynomial chaos expansion.
The $\mathbb{P}^N$ represents the polynomial space spanned by $\{\phi_n(\xi)\}$ for $n=0,\, ... \, , \,N$.

Importantly, the free boundary $S(t,\xi)$ is not explicitly available from the mathematical model but must be numerically obtained by solving the enthalpy formulation. To examine the impact of uncertain parameters on the free boundary, we need to quantify how the randomness in $\xi$ propagates through the model to $S(t,\xi)$.

To achieve this, we leverage a generalized polynomial chaos (gPC) approach (cf. \cite{book-Stefan-UQ,UQflight}). The key idea is to represent $S(t,\xi)$ as an expansion in terms of polynomial basis functions $\phi_n(\xi)$:

  	\begin{align}\label{eq:UQ1}
  	S(t,\xi)=\sum_{n=0}^{\infty}c_n(t)\phi_n(\xi),
  \end{align}
  	where $c_n(t)$ are the gPC coefficients and the index of summation $n$ indicates the polynomial order.  This provides a surrogate model relating the random inputs $\xi$ to the quantity of interest $S(t,\xi)$. We can then analyze this surrogate to quantify uncertainty propagation through the model.
   The selection of the polynomial basis function depends on the distribution of 
   the uncertainty parameters (see, cf. \cite{book-Stefan-UQ}).
  In practical use, we usually truncate it with finite $N$ terms. 
  	\par
   Different approaches are available to determine the gPC coefficients $c_n$. We opt to leverage a non-intrusive stochastic collocation approach rather than an intrusive formulation such as the stochastic Galerkin method \cite{UQGalerkin1,UQGalerkin2,UQGalerkin3}. Though intrusive methods can directly embed uncertainty into the model, they require extensive modifications to the enthalpy equations and discretization for this Stefan problem. The presence of two moving boundaries (injection and phase interface) evolving simultaneously poses challenges for choosing appropriate polynomial chaos basis functions when uncertainties are embedded, due to discontinuities and sharp solution gradients. The complexity of this problem imposes significant additional hurdles for an intrusive approach.

In contrast, the non-intrusive stochastic collocation method only requires performing sampling and simulations using the original enthalpy solver, without needing to modify it. With a sufficient number of solution observations, the gPC expansion coefficients can be determined by post-processing the simulation data. While both methods face potential convergence difficulties, the implementation challenges of intrusive UQ for this  Stefan problem make a non-intrusive approach leveraging the existing enthalpy solver more attractive. The modular nature avoids extensive changes to the solver while still allowing uncertainty propagation.

    In this work, from a comprehensive consideration of simplicity and effectiveness,
 we employ the stochastic collocation method and apply the discrete Least Squares projection approach to determine
coefficients of (\ref{eq:UQ1}).
 For each $\xi_m$, we discretize time $t$ and then obtain $(S(t_1, \, \xi_m),\, \cdots \, , \,S(t_k, \, \xi_m))$.
Then we can  determine the coefficients $c_n$ by solving the following optimization problem which minimizes the approximation error:
	\begin{align}\label{eq:UQ2}
		S_N(t_k,\xi):=P_M^N S=\operatorname{arg}\min\limits_{u\in\mathbb{P}^N}\frac{1}{M}\sum\limits_{m=1}^M\left(u(t_k,\xi_m)-S(t_k,\xi_m)\right)^2,~ k=1,..,K.
	\end{align}
    We introduce the discrete inner product as
    \begin{align}
    	\langle u,v\rangle_M=\frac{1}{M} \sum_{m=1}^{M} u(\xi_m)v(\xi_m).
    \end{align}
    Then we can rewrite equation (\ref{eq:UQ2}) as
    \begin{align}
   S_N(t_k,\xi)=\underset{u\in\mathbb{P}^{N}}{\operatorname{arg}\min}\|u-S\|_M,
   \end{align}
    where $\|u\|_M=\langle u,u\rangle_M^{1/2}.$ 
    It is equivalent to solve 
    \begin{align} 
    	\hat{A}c_k=\hat{S_k}, \quad k=1,...K.
	\end{align}
	where 
	\begin{align}
		\hat{A}=(\langle \phi_i,\phi_j\rangle)_{N \times N},\quad \hat{S_k}=(\langle S(t_k),\phi_j\rangle)_{N \times 1} ,\quad k=1,..K.
	\end{align}

     As we mentioned above, the stochastic collocation method employs sample information to determine unknown coefficients ${c_n}$ with a high degree of precision, which only requires solving multiple equations of the same size and form as the original model.
     There are also other methods to determine
$\hat{A}$ and $\hat{S(t_k)}$.
     Further details can be found in \cite{UQcollocation1}.\

  The uncertainty quantification techniques presented thus far focused on the 1D mathematical model and computational framework. As formulated in Sec. \ref{sec:model}, the model naturally extends to 2D scenarios on a domain $(y,x) \in [0,1] \times [0,L(t)]$.
In the 2D setting, the free boundary becomes a curve dependent on both spatial dimensions and time, denoted $S(y,t,\xi')$. Note that $\xi'=f(y,\xi)$ may have $y$ dependence. The evolution of $S(y,t,\xi')$ depends on uncertainty in parameters $\xi$ as well as the additional spatial variable $y$. Then equation (\ref{eq:UQ1}) can be changed to:
\begin{align}
S(y,t,\xi')=\sum_{n=0}^{\infty}c_n(t,y)\phi_n(\xi).
\end{align}

For each sampled $\xi_m$, we first discretize the spatial variable $y$ into a set of points $\{y_j\}_{j=1}^{J}$. Next, we discretize time into a set of points $\{t_k\}_{k=1}^{K}$. Numerical solution of the enthalpy system at each combination of $y_j$ and $t_k$ yields the free boundary values $(S(y_j,\, t_1, \, \xi_m'),\, \cdots \, , \, S(y_j,\, t_K,\, \xi_m'))$ for the given $\xi_m'$. Applying polynomial chaos regression at each spatial point $y_j$ determines the gPC coefficients $\{c_n(t_k,y_j)\}_k$ analogous to the 1D case. Repeating this process for each $y_j$ provides the full set of gPC coefficients $\{c_n(t_k,y_j)\}_{k,j}$ describing the surrogate model for the 2D free boundary $S^{\star}(y,t,\xi')$.

The overall UQ methodologies remain applicable, but the increase in physics dimensionality and geometric complexity introduces new demands on sampling, approximations, analysis, and computation. More samples may be needed to sufficiently cover the higher dimensional space. However, the general UQ framework provides a starting point that can be built upon with extensions addressing these 2D considerations. 

In practice, we need to determine the number of samples $M$ and the truncated expansion order $N$ together. It is shown in \cite{UQuniform} that when samples are randomly generated from the probability density function, choosing sample size $M$ scaling quadratically with $N$ (i.e. $M=cN^2$) leads to a stable discrete least squares projection with high probability. This provides practical guidance on selecting suitable values of $M$ and $N$ for the gPC approximation to achieve stable UQ results.
	
	Putting together the procedures above, we can then present the whole 1D process as it shown in Algorithm \ref{algor:UQ}. In 2D settings, we only need to repeat steps 4 to 11  of Algorithm \ref{algor:UQ}  $J$ times to get $S$ for all $t_k$ and all $y_j$.
	\vskip 3mm
 \begin{algorithm}
\caption{Uncertainty Quantification Approach for 1D}

\begin{algorithmic}[1]

\STATE Choose polynomial chaos basis ${\phi_i(\xi)}$.

\STATE Determine number of samples $M$.
\STATE Generate input samples $\{\xi^m\}_{m=1}^M$.

\FOR{$m = 1$ \TO $M$}
\STATE Sample parameters: $\xi \leftarrow \xi^m$.
\STATE Run simulator to obtain $S(t_k, \xi^m)$.
\ENDFOR

\STATE Compute inner products: $\hat{A} = (\langle\phi_i,\phi_j\rangle)$, $\hat{S}_k = (\langle S(t_k),\phi_j\rangle)$.
\FOR{each $t_k$}
\STATE Solve linear system $\hat{A}c_k = \hat{S}_k$.
\ENDFOR

\STATE Construct surrogate: $S_N(t, \xi) = \sum_i c_i(t) \phi_i(\xi)$.

\STATE Compute statistics from surrogate model

\end{algorithmic}
\label{algor:UQ}
\end{algorithm}
\section{Numerical Tests and Experiments}\label{sec:Numer}
In this section, we present both one and two-dimensional numerical examples to verify the effectiveness of the proposed models and algorithms and to explore the  various modeling phenomena, such as the determination of icing type, the effect of horizontal thermal diffusion, and the IHI parameter on ice formation shape,
and the influence of uncertainty quantification with random parameters, etc. Note that our discussion will focus on the particular scenario of airfoil icing, but other interpretations are possible when other phase transition scenarios are considered.

\subsection{Numerical Experiments in One-dimension}\label{subsec:5.1}
The model as specified in Eq. (\ref{model:1}) is considered first. The properties and the conditions considered were:
 \begin{equation}
\begin{aligned}
    &\rho=1 kg/m^3, ~ K=2W/m \cdot ^\circ \mathrm{C}, ~
  C=2.5 MJ/kg \cdot ^\circ \mathrm{C} ,~ L_h=100 MJ/kg,~ \\
  &T_{\text{\text{initial}}}=
  2^\circ \mathrm{C},~ T_{0}=-10^\circ \mathrm{C},~ 
  T_{m}=0^\circ \mathrm{C},~ L_{0}=1 m, ~L_{\text{ref}}=1 m.
\end{aligned}
 \end{equation}

\begin{figure}[ht] 

  \centering
  \subfigure[\( \widehat{\eta}=1.25\), \(\widehat{\beta}=0.35\), \(\left(S^{\star}(4)\approx 1.15\right) \)]{
  \begin{minipage}[t]{0.47\textwidth}
  \centering
  \includegraphics[width=6cm]{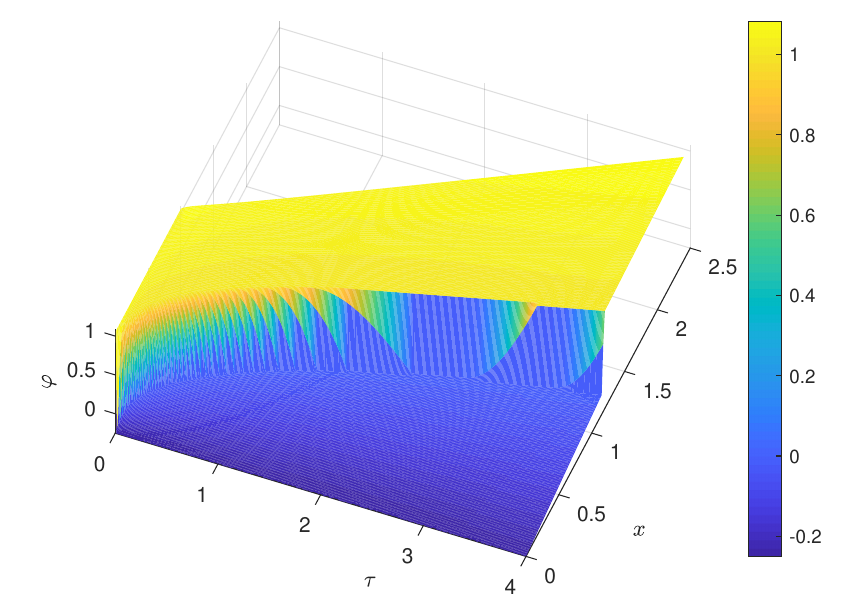}
  \label{fig:eta(a)}
  \end{minipage}
  }
  \subfigure[\( \widehat{\eta}=1.25\), \(\widehat{\beta}=1 \), \(\left(S^{\star}(4)\approx 1.05\right) \)]{
  \begin{minipage}[t]{0.47\textwidth}
  \centering
  \includegraphics[width=6cm]{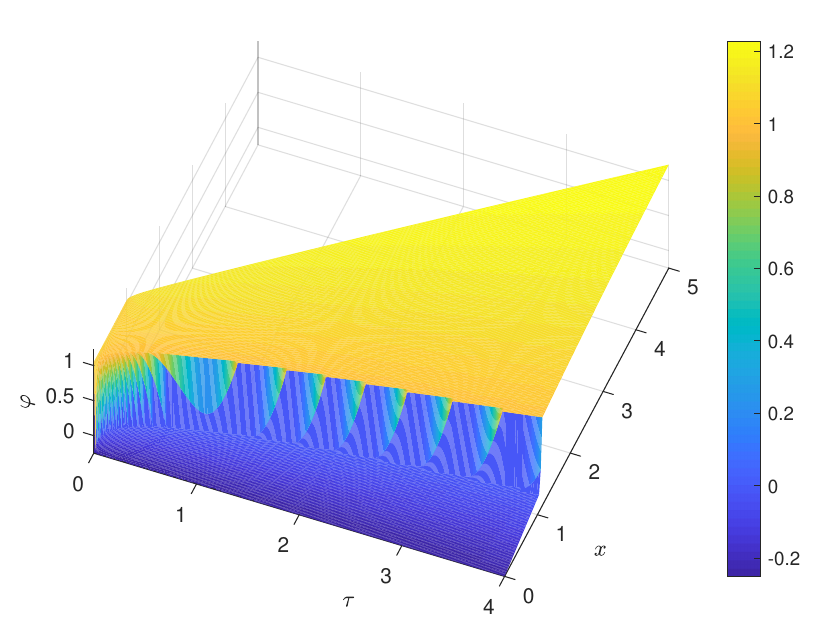}
  \label{fig:eta(b)}
  \end{minipage}
  }
  \subfigure[\( \widehat{\eta}=1\), \(\widehat{\beta}=0.35 \), \(\left(S^{\star}(4) \approx 1.32\right) \)]{
  \begin{minipage}[t]{0.47\textwidth}
    \centering
    \includegraphics[width=6cm]{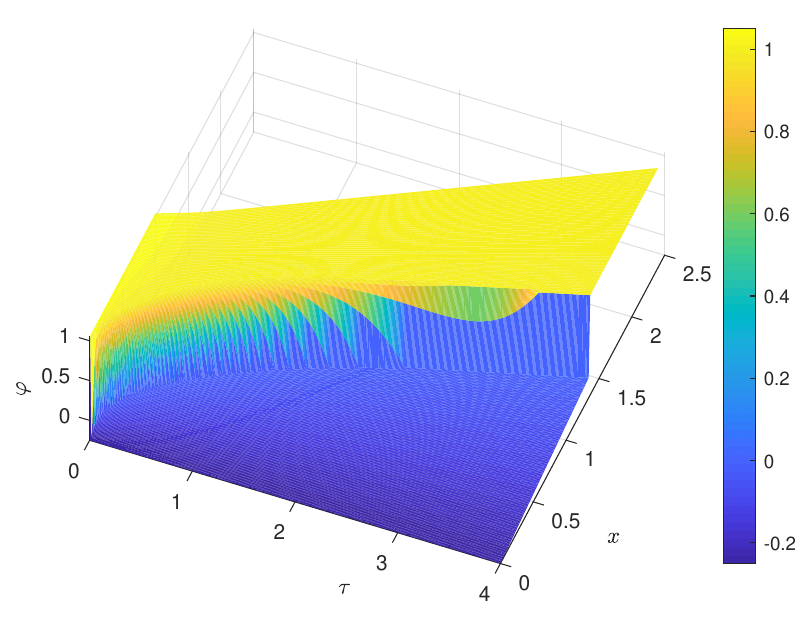}
    \end{minipage}
  }
  \label{fig:EX1D}
  \subfigure[\( \widehat{\eta}=1\), \(\widehat{\beta}=1 \), \(\left(S^{\star}(4) \approx 1.30\right) \)]{
    \begin{minipage}[t]{0.47\textwidth}
      \centering
      \includegraphics[width=5.7cm]{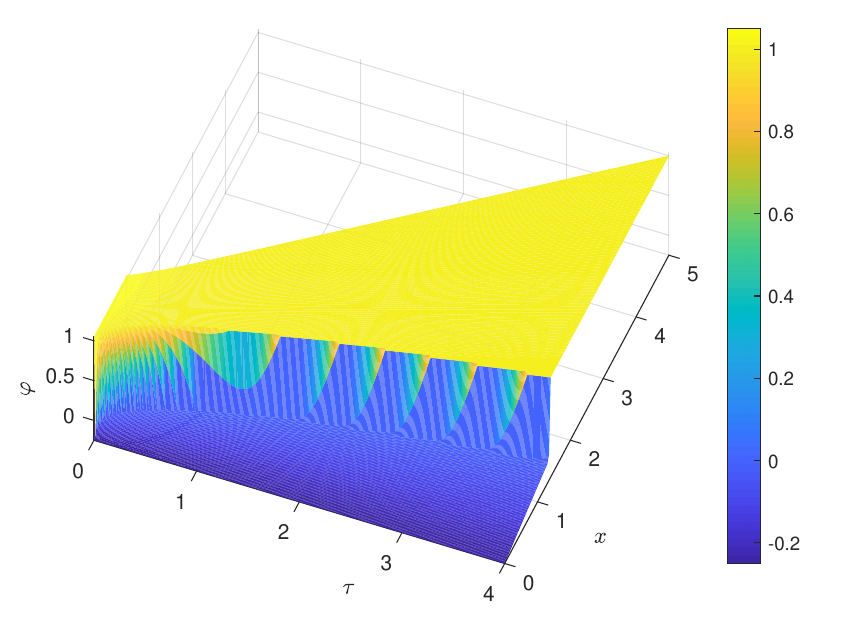}
      \end{minipage}
  }
      \caption{  Example in One-dimension: The enthalpy value with the 
      mesh size \(\Delta x\) being 0.01 and \(\Delta \tau\) 
      being 0.0001 and time $\tau$
      changing from 0 to 4.  The free boundary location $S^{\star}(\tau)$ at $\tau=4$.}

\end{figure}

The Fig. \ref{fig:EX1D} shows the evolution of the enthalpy profile over time in the one-dimensional system.
The horizontal axis is time  $\tau$ ranging from 0 to 4, as the injection boundary expands with time.  The vertical axis shows the enthalpy value.
The spatial discretization used $\Delta x = 0.01$ and time step $\Delta \tau = 0.0001$ with different choices of  incoming heat influx 
(IHI)   \(\widehat{\eta}\) and injection boundary velocity
(IBV)   \(\widehat{\beta}\).

At any given $\tau$ snapshot, the enthalpy exhibits a discontinuity at the inner interface between solid and liquid phases. Tracking this jump point visualizes the competing growth of the freezing interface and injection boundary.

Several key features manifest through the simulation results.
Initially, the inner freezing interface exhibits rapid growth, surpassing the velocity of the injection boundary. This indicates the freezing rate is outpacing the water influx rate. We can interpret this phase of icing as a physical scenario where incoming water droplets freeze instantaneously upon contact with the solid surface. Such behavior is characteristic of rime ice formation. Note that the link between the interface behavior and these terms is rather phenomenological, since in aircraft icing, the physical process is extremely complicated and cannot be fully recovered by our model. 
At later times, the inner interface velocity slows down while the injection boundary continues expanding linearly. This overflow behavior might be indicative of the so-called glaze ice formation, typically occurring at warmer temperatures. Here, the water influx rate starts exceeding the freezing rate, leading to liquid water overflowing past the existing ice layer.

In summary, the numerical experiments presented demonstrate the model's efficacy in capturing complex morphological physics underlying rime and glaze ice formation. The time evolution of the inner interface and outer boundary locations, along with their relative velocities, agrees with established theory regarding the competing effects of thermal diffusion and mass influx. The results showcase the model’s ability to reproduce transitions from complete freezing to water accumulation regimes under various conditions. In particular, the dependence on parameters such as outer boundary velocity and incoming droplet energy content produces quantitative and qualitative outcomes consistent with physical intuition.

    \subsection{Numerical Experiments in Two-dimension}\label{subsec:Num_2D}
Now we consider the Eq. (\ref{eq:category}) in two-dimension.
In the one-dimensional numerical example, we study the influence of parameters on
 the interface moving speed,  so in the two-dimensional numerical example, we are more like to study the influence of parameters and the effect of horizontal thermal diffusion on the ice formation shape. Some conditions and parameters in this example were:

  \begin{equation}
    \begin{aligned}
      &\rho=1 kg/m^3, ~ K=2W/m \cdot ^\circ \mathrm{C}, ~
    C=2.5 MJ/kg \cdot ^\circ \mathrm{C}, ~ L_h=100 MJ/kg, ~ \\
    &T_{\text{\text{initial}}}=
    4^\circ \mathrm{C}, ~ T_{0}=-4^\circ \mathrm{C}, ~
    T_{m}=0^\circ \mathrm{C}, ~ L_{0}=0.2 m, ~ \widehat{\beta}=0.1, ~ L_{\text{ref}}=1 m.
  \end{aligned}
  \end{equation}

\begin{figure}[ht]
      
  \centering
  \subfigure[]{
  \begin{minipage}[t]{0.45\textwidth}
  \centering
  \includegraphics[width=6cm]{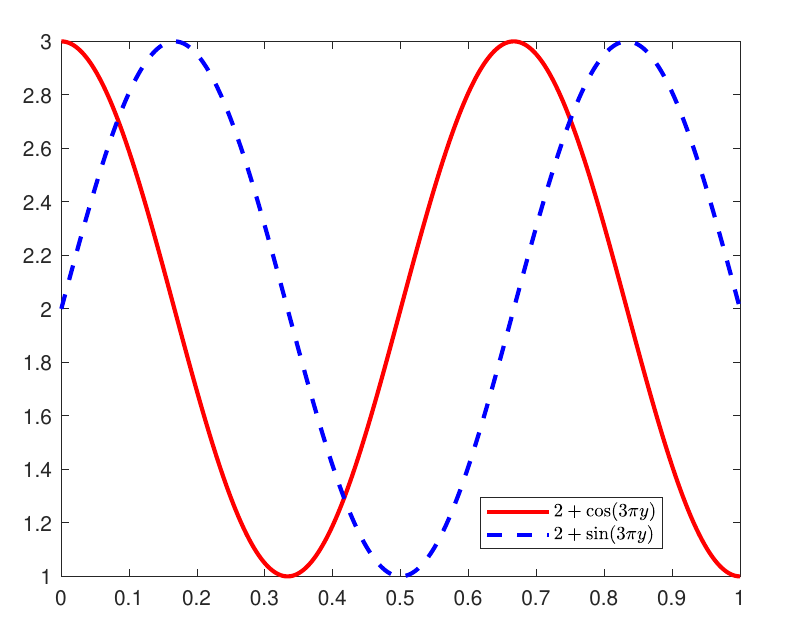}
  \label{fig:eta_func1}
  \end{minipage}
  }
  \subfigure[]{
  \begin{minipage}[t]{0.45\textwidth}
  \centering
  \includegraphics[width=6cm]{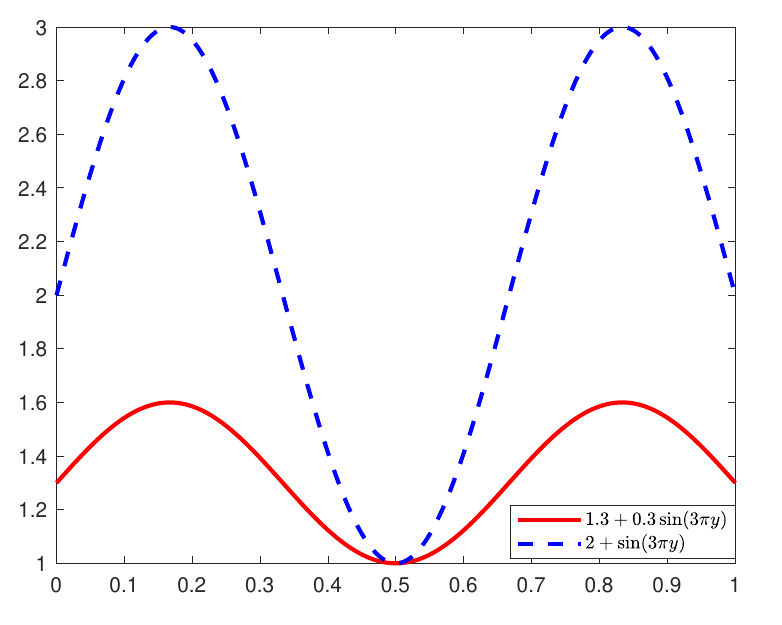}
  \label{fig:eta_func2}
  \end{minipage}
  }
      \caption{Example in Two-dimension: The different choice of 
      \(\widehat{\eta}\).}
\end{figure}

In the two-dimensional case, the IHI parameter \(\widehat{\eta}\) may not be homogeneous in the horizontal \(y\) direction. So we study the effects of different choices of IHI parameter \(\widehat{\eta}\)  on the shape of ice formation, such as the different form, as shown in Fig. \ref{fig:eta_func1} or the same form, but different amplitudes, as shown in Fig. \ref{fig:eta_func2}. 
We take $\widehat{\eta}=2+\cos(3\pi y)$,  $2+\sin(3\pi y)$, and $2+\sin(3\pi y)$ respectively.

\begin{figure}[ht]
  \centering
  \subfigure[\(\tau=0.3\)]{
  \begin{minipage}[t]{0.45\textwidth}
  \centering
  \includegraphics[width=6cm]{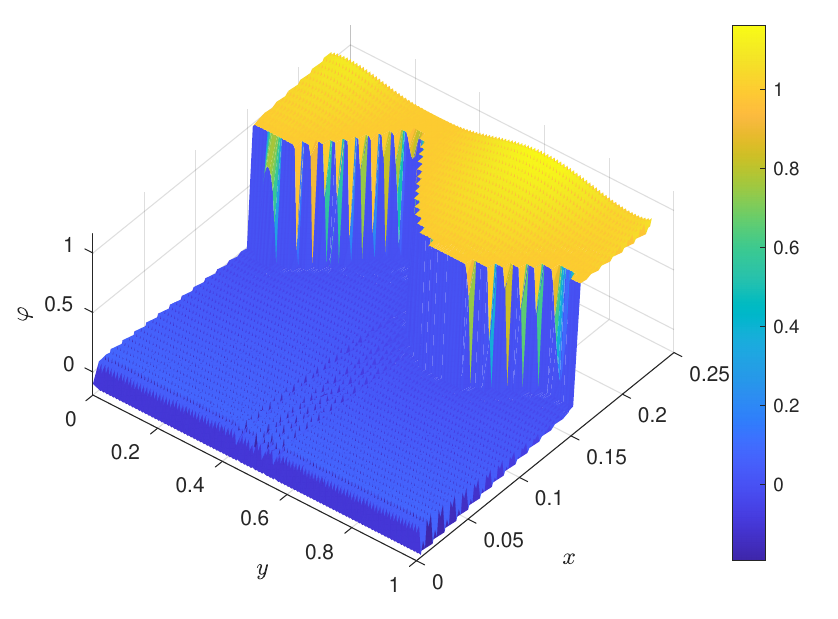}
  \label{Fig:eta1(a)}
  \end{minipage}
  }
  \subfigure[\(\tau=0.7\)]{
  \begin{minipage}[t]{0.45\textwidth}
  \centering
  \includegraphics[width=6cm]{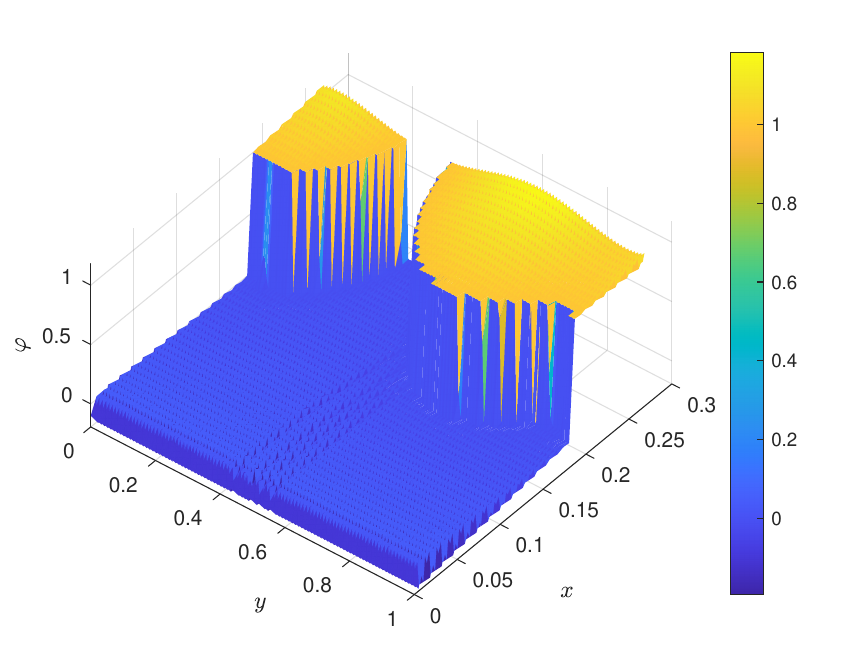}
  \label{Fig:eta1(b)}
  \end{minipage}
  }
  \subfigure[\(\tau=1\)]{
  \begin{minipage}[t]{0.45\textwidth}
    \centering
    \includegraphics[width=6cm]{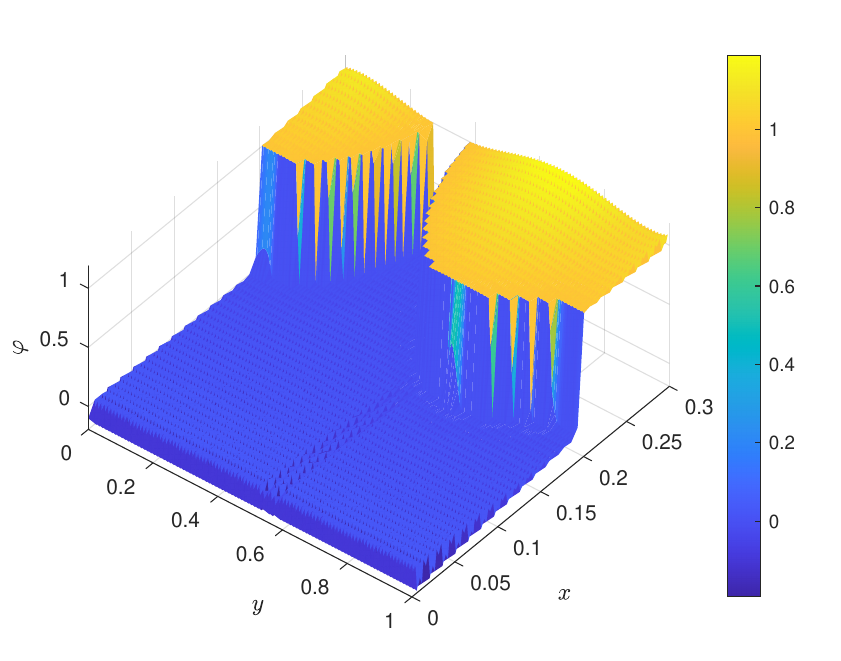}
  \label{Fig:eta1(c)}
    \end{minipage}
  }
  \subfigure[\(\tau=2\)]{
    \begin{minipage}[t]{0.45\textwidth}
      \centering
      \includegraphics[width=5.7cm]{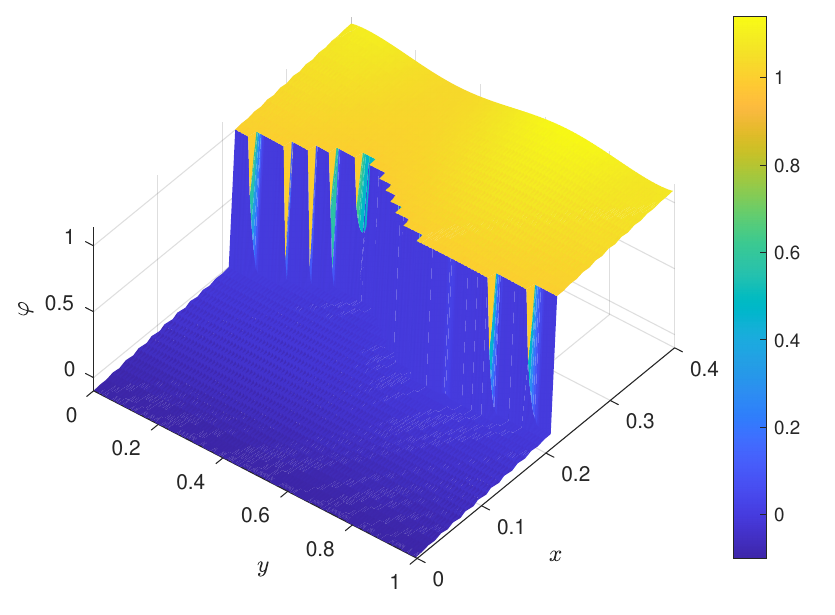}
  \label{Fig:eta1(d)}
      \end{minipage}
  }
      \caption{Example in Two-dimension: 
      The enthalpy value and the icing shape with 
      \(\widehat{\eta}=2+\cos(3\pi y)\),
      the mesh size \(\Delta y\), \(\Delta z\) being 0.005, the 
      time step \(\Delta \tau\) 
      being 0.00001 and the time taking 0.3, 0.7, 1, 2 
      respectively. The blue area corresponds to  liquid while 
      the yellow area corresponds to solid.}
      \label{Fig:eta1}
\end{figure}

\begin{figure}[ht]
      
  \centering
  \subfigure[\(\tau=0.3\)]{
  \begin{minipage}[t]{0.45\textwidth}
  \centering
  \includegraphics[width=6cm]{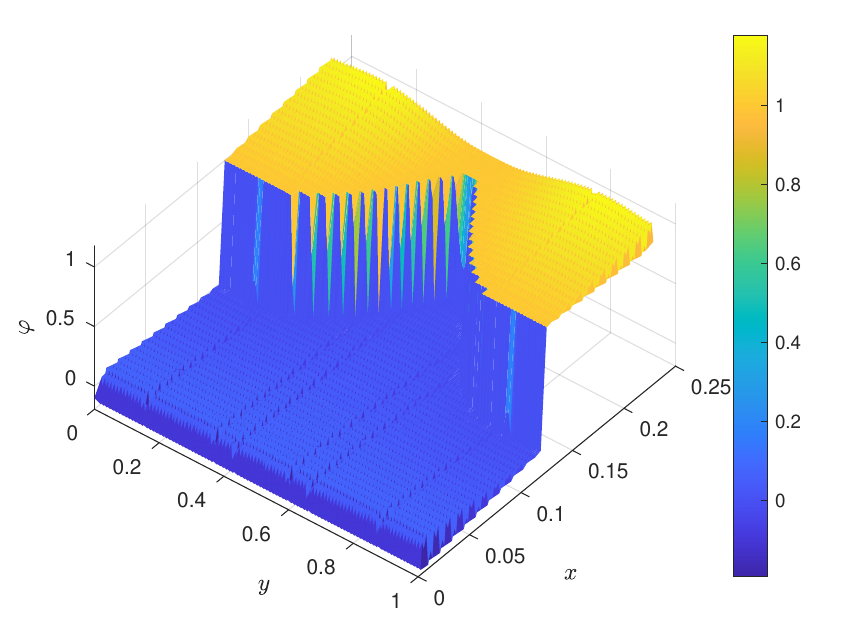}
  \end{minipage}
  }
  \subfigure[\(\tau=0.7\)]{
  \begin{minipage}[t]{0.45\textwidth}
  \centering
  \includegraphics[width=6cm]{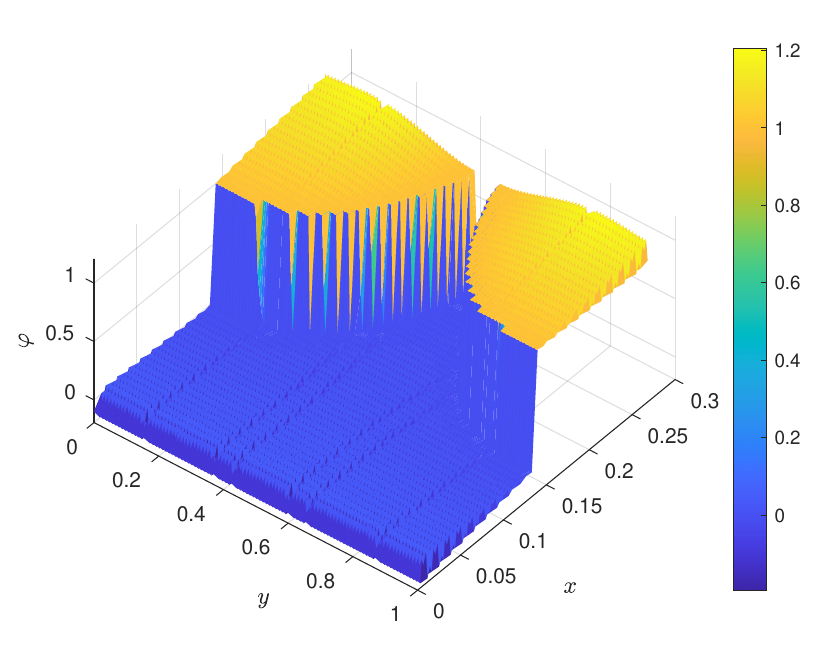}
  \end{minipage}
  }
  \subfigure[\(\tau=1\)]{
  \begin{minipage}[t]{0.45\textwidth}
    \centering
    \includegraphics[width=6cm]{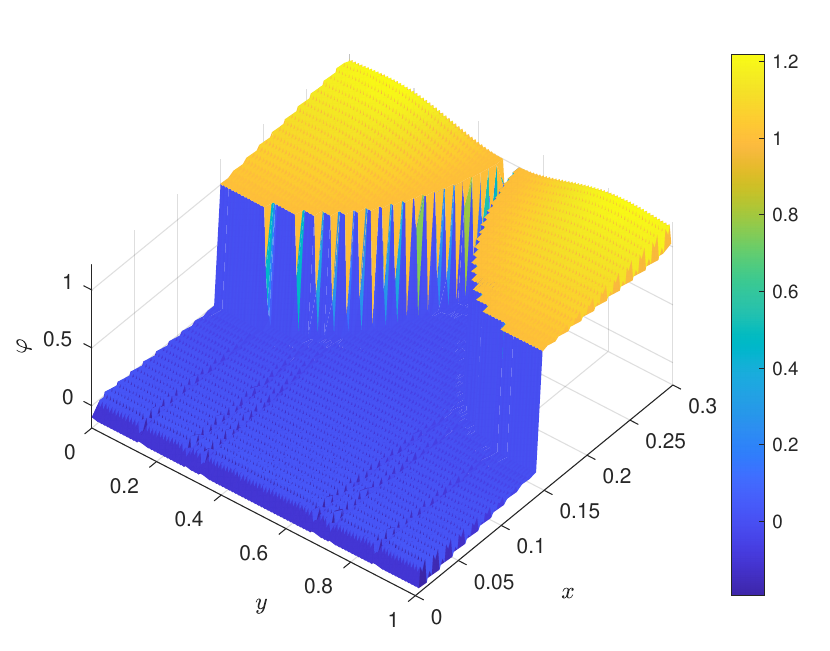}
    \end{minipage}
  }
  \subfigure[\(\tau=2\)]{
    \begin{minipage}[t]{0.45\textwidth}
      \centering
      \includegraphics[width=5.7cm]{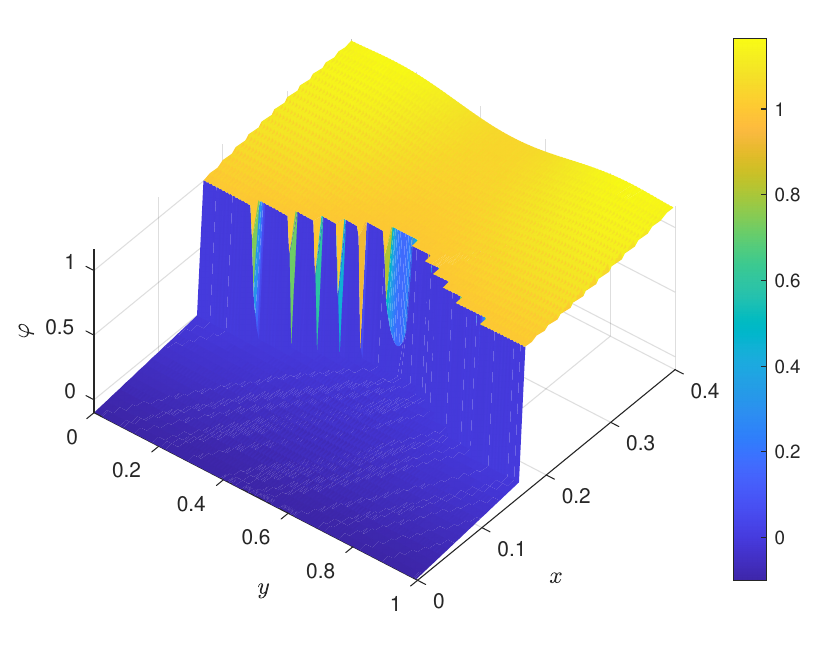}
      \end{minipage}
  }
      \caption{Example in Two-dimension: 
      The enthalpy value and icing shape with 
      \(\widehat{\eta}=2+\sin(3\pi y)\),
      the mesh size \(\Delta y\), \(\Delta z\) being 0.005, the 
      time step \(\Delta \tau\) 
      being 0.00001 and the time taking 0.3, 0.7, 1, 2 
      respectively. The blue area corresponds to  liquid while 
      the yellow area corresponds to solid.}
      \label{fig:eta2}
\end{figure}

\begin{figure}[ht]
      
  \centering
  \subfigure[\(\tau=0.1\)]{
  \begin{minipage}[t]{0.45\textwidth}
  \centering
  \includegraphics[width=6cm]{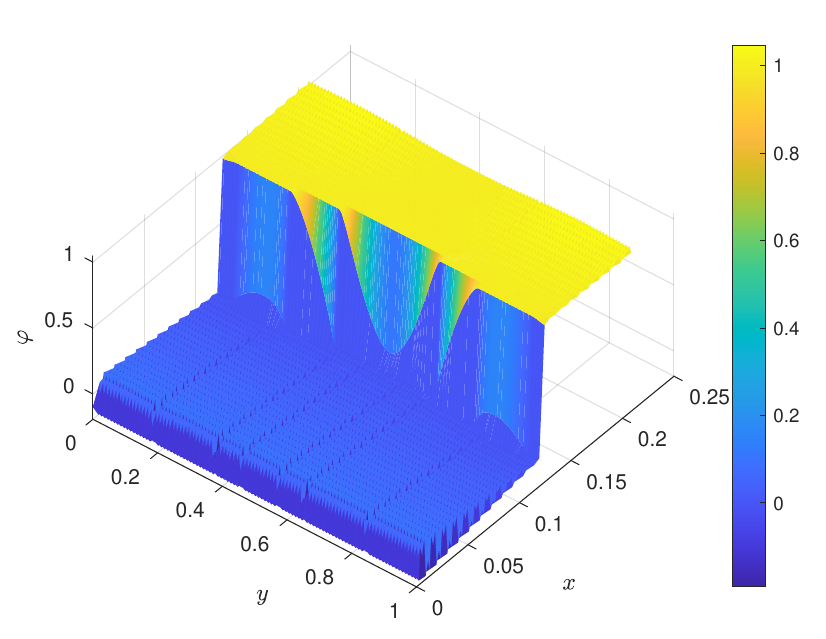}
  \end{minipage}
  }
  \subfigure[\(\tau=0.2\)]{
  \begin{minipage}[t]{0.45\textwidth}
  \centering
  \includegraphics[width=6cm]{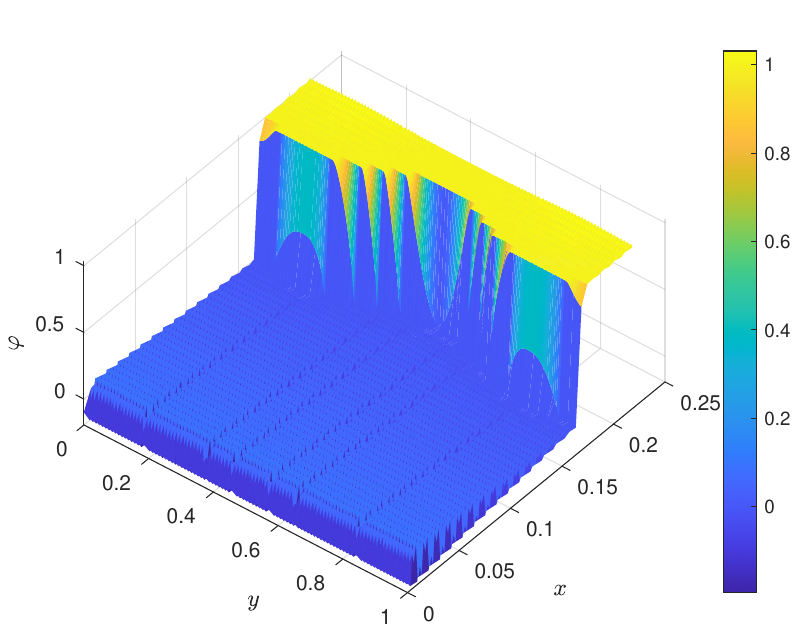}
  \end{minipage}
  }
  \subfigure[\(\tau=0.3\)]{
  \begin{minipage}[t]{0.45\textwidth}
    \centering
    \includegraphics[width=6cm]{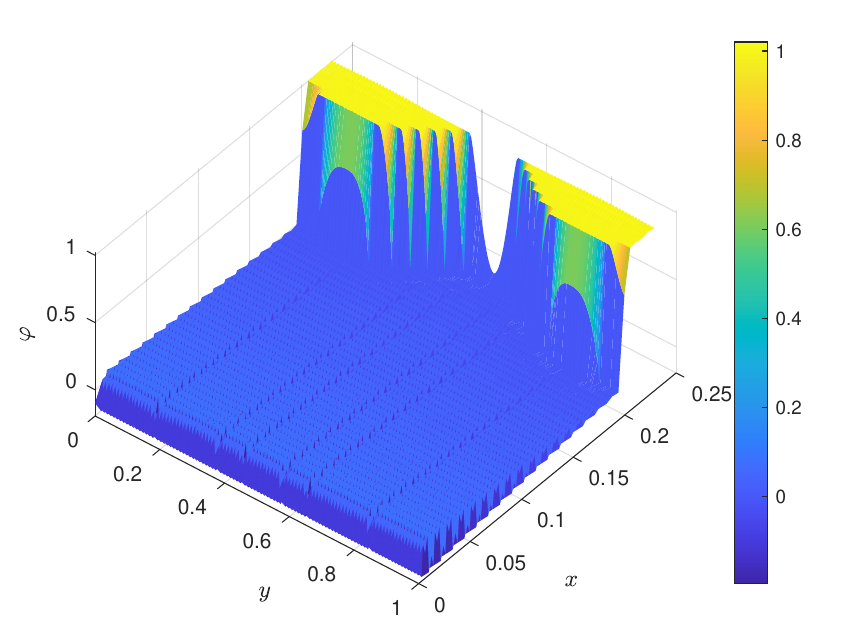}
    \end{minipage}
  }
  \subfigure[\(\tau=0.5\)]{
    \begin{minipage}[t]{0.45\textwidth}
      \centering
      \includegraphics[width=5.7cm]{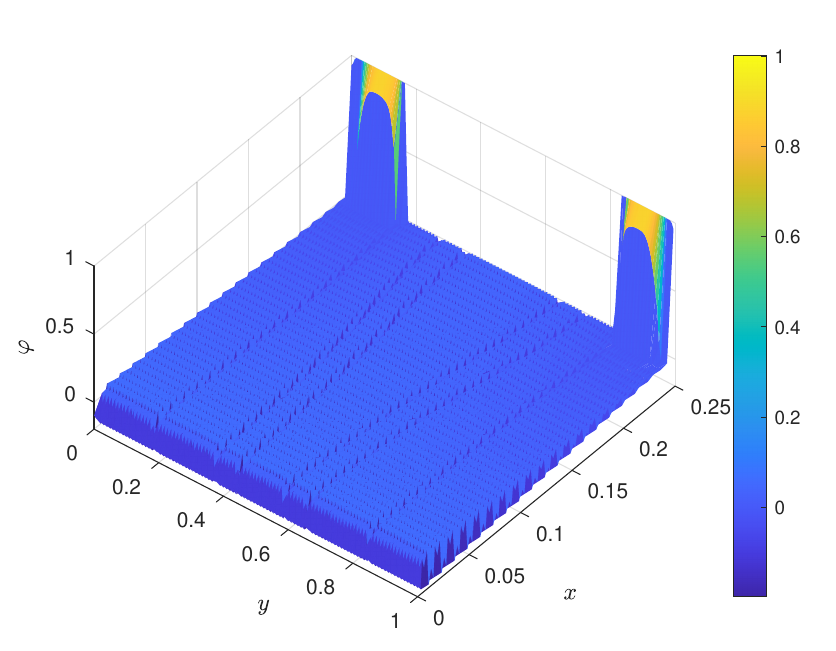}
      \end{minipage}
  }
      \caption{Example in Two-dimension: 
      The enthalpy value and icing shape with 
      \(\widehat{\eta}=1.3+0.3\sin(3\pi y)\),
      the mesh size \(\Delta y\), \(\Delta z\) being 0.005, the 
      time step \(\Delta \tau\) 
      being 0.00001 and the time taking 0.1, 0.2, 0.3, 0.5 
      respectively. The blue area corresponds to  liquid while 
      the yellow area corresponds to solid.}
      \label{fig:eta3}
\end{figure}

The Fig. \ref{Fig:eta1}-\ref{fig:eta3} depict the simulated evolution of ice morphology on a 2D airfoil using an enthalpy-based model. 
The $y$ and $z$ axes represent the spatial coordinates
and the vertical axis shows the enthalpy value.
The ice interface is delineated by a color contour, with blue and yellow indicating liquid water and fully frozen ice respectively. This visualizes the phase boundary location over time.
The model was discretized with $\Delta y=\Delta z=0.005$ and time step $\Delta \tau=0.00001$ with different choices of incoming water droplet IHI parameter $\widehat{\eta}$. The subfigures show the simulated ice shape at different dimensionless times. 

We observe that since the IHI parameter \(\widehat{\eta}\) is not 
homogeneous, the ice formation interface is not homogeneous either, and the corresponding shape is consistent with the form of \(\widehat{\eta}\).
From our numerical results,  we can see that ice formation can be roughly divided into three stages.

Initially in Fig. \ref{Fig:eta1(a)}-\ref{Fig:eta1(b)}, the interface is close to the boundary, inducing rapid freezing that overcomes even high energy droplet areas, characteristic of rime ice. Subsequently in Fig. \ref{Fig:eta1(b)}-\ref{Fig:eta1(c)}, the weakening thermal gradient allows local energy variations to manifest. High energy areas begin resisting freezing, transitioning toward glaze ice formation as latent heat exceeds reducing heat transfer. Finally in Fig. \ref{Fig:eta1(c)}-\ref{Fig:eta1(d)}, the phase change interface moves further from the boundary, reducing its freezing impact. The incoming droplet energy now exceeds the freezing driven by the airfoil boundary, evidenced by the decreasing amplitude of the heterogeneous ice shape over time as surface effects diminish. The system transitions toward a latent heat-dominated state.

   In Fig. \ref{Fig:eta1}-\ref{fig:eta2}, we can see that different forms of IHI parameters lead to different ice formation shapes, but the ice formation process meets our three stages of ice formation described above.

   In Fig. \ref{fig:eta2}-\ref{fig:eta3},  the \(\widehat{\eta}\) is taken in the form of a \(\sin\) function, while the amplitude is taken differently. 
   We found that in the smaller amplitude case (Fig. \ref{fig:eta3}), due to the small amount of incoming water droplet energy, the impact brought by it is significantly smaller than the impact brought by the airfoil boundary. 

    The two-dimensional numerical experiments validate the model’s capacity to capture complex icing physics using the spatially- and temporally-varying latent heat parameter $\widehat{\eta}(y,t)$. Specifically, the results showcase the model’s capacity to reproduce heterogeneous ice shapes arising from non-uniform energetic and transport conditions. As in the 1D case, the quantitative evolution matches established qualitative expectations. Overall, the preliminary 2D results provide promising indications that this modeling approach can provide new insight into the nuanced physics underlying the multifaceted morphology of atmospheric ice accretion. And ongoing efforts should focus on more rigorous verification and validation using experimental data.

\subsection{Uncertainty Quantification in One-dimension}\label{subsec:UQ}
Consider the uncertainty quantification process described in Sec. \ref{sec:UQ}. 
 The physical properties and conditions remain the same as Sec. \ref{subsec:5.1}, except for the airfoil surface temperature $T_{0}$, incoming boundary velocity (IBV)  $\widehat{\beta}$, and incoming heat influx (IHI)  $\widehat{\eta}$. 
 The IBV and IHI parameters are treated as independent uniform random variables to quantify their uncertainty impacts on the phase change interface location.
 The finite difference mesh size is $\Delta x = 0.01$ and the time step is $\Delta \tau = 0.0001$, with the time $\tau$ spanning 0 to 4.
Fourth-order Legendre polynomials are utilized as the gPC basis functions to construct surrogate models for the interface location.

\begin{figure}[ht]
      
  \centering
  \subfigure[$T_{0}=-26^\circ C$, $\mathbb{E}\left(S^{\star}(4)\right)\approx1.949$, $\sigma\approx0.067$]{
  \begin{minipage}[t]{0.45\textwidth}
  \centering
  \includegraphics[width=6cm]{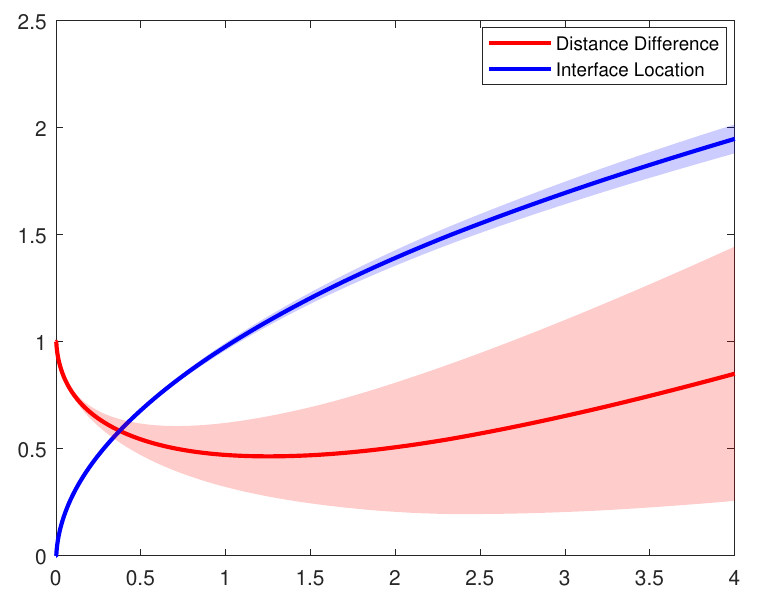}
  \end{minipage}
  }
  \subfigure[$T_{0}=-18 ^\circ C$, $\mathbb{E}\left(S^{\star}(4)\right)\approx 1.651$, $\sigma\approx 0.066$]{
  \begin{minipage}[t]{0.45\textwidth}
  \centering
  \includegraphics[width=6cm]{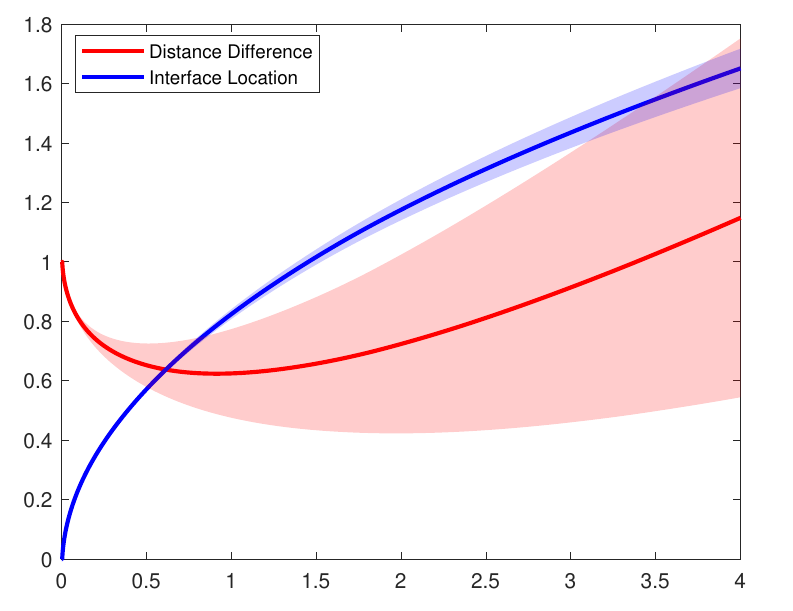}
  \end{minipage}
  }
  \subfigure[$T_{0}=-10 ^\circ C$, $\mathbb{E}\left(S^{\star}(4)\right)\approx1.239$, $\sigma\approx0.061$]{
  \begin{minipage}[t]{0.45\textwidth}
    \centering
    \includegraphics[width=6cm]{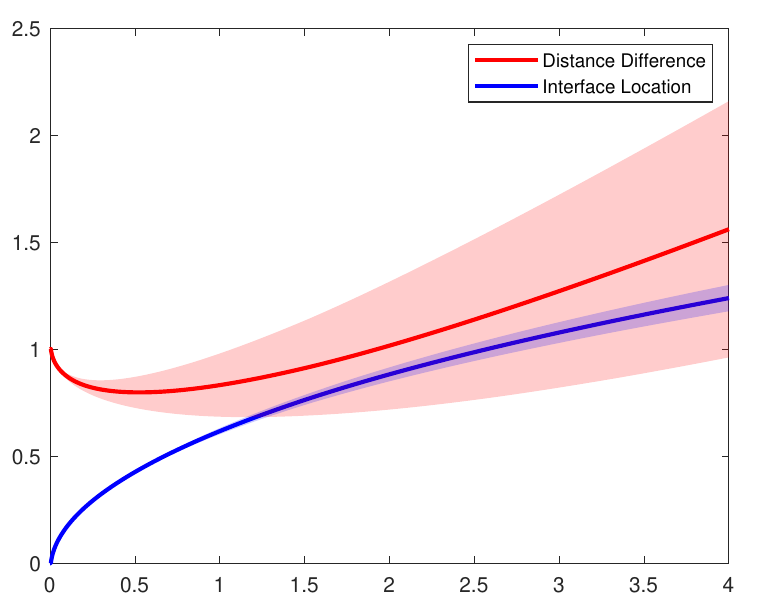}
    \end{minipage}
  }
  \subfigure[$T_{0}=-2 ^\circ C$, $\mathbb{E}\left(S^{\star}(4)\right)\approx0.520$, $\sigma\approx0.049$]{
    \begin{minipage}[t]{0.45\textwidth}
      \centering
      \includegraphics[width=5.7cm]{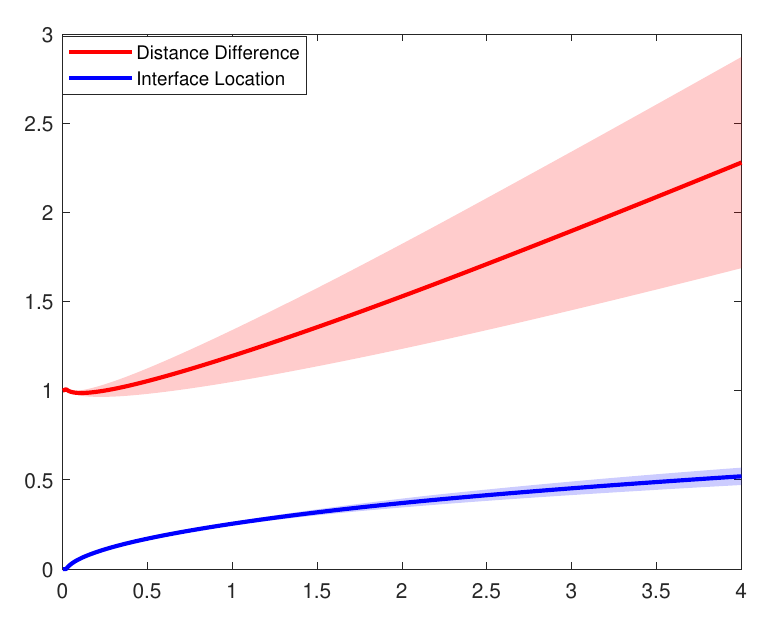}
      \end{minipage}
      
  }
\caption{UQ in One-dimension: The phase change interface location and distance between inner and outer boundaries are shown for influx uncertainty modeled by $\widehat{\beta} \sim \mathcal{U}(0.2,0.7)$ and $\widehat{\eta} \sim \mathcal{U}(1,1.25)$ and different $T_0$ options. Shaded regions indicate plus/minus one standard deviation uncertainty bands on the interface location and boundary separation distance. The mean and standard deviation of the free boundary location $S^{\star}(\tau)$ at $\tau=4$ are also presented.}
        \label{wingtemp}
  \end{figure}

\begin{figure}[ht]
      
  \centering
  \subfigure[$\widehat{\eta} \sim \mathcal{U}(1,1.25)$, $\mathbb{E}\left(S^{\star}(4)\right)\approx 1.239$, $\sigma\approx 0.060$]{
  \begin{minipage}[t]{0.45\textwidth}
  \centering
  \includegraphics[width=6cm]{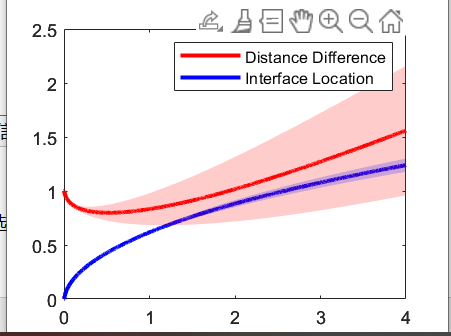}
  \end{minipage}
  }
  \subfigure[$\widehat{\eta} \sim \mathcal{U}(1.25,1.5)$,  $\mathbb{E}\left(S^{\star}(4)\right)\approx 1.061$, $\sigma \approx 0.065$]{
  \begin{minipage}[t]{0.45\textwidth}
  \centering
  \includegraphics[width=6cm]{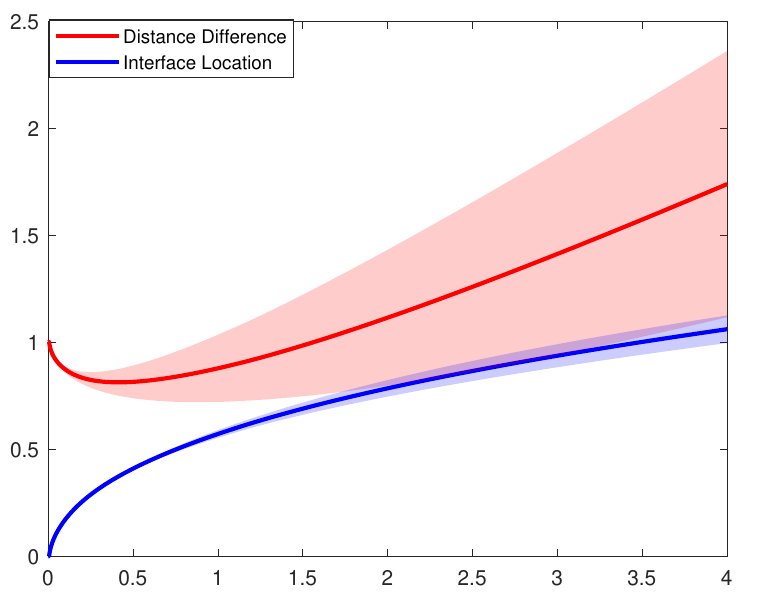}
  \end{minipage}
  }
  \subfigure[$\widehat{\eta} \sim \mathcal{U}(1.5,1.75)$,  $\mathbb{E}\left(S^{\star}(4)\right)\approx 0.917$, $\sigma\approx0.072$]{
  \begin{minipage}[t]{0.45\textwidth}
    \centering
    \includegraphics[width=6cm]{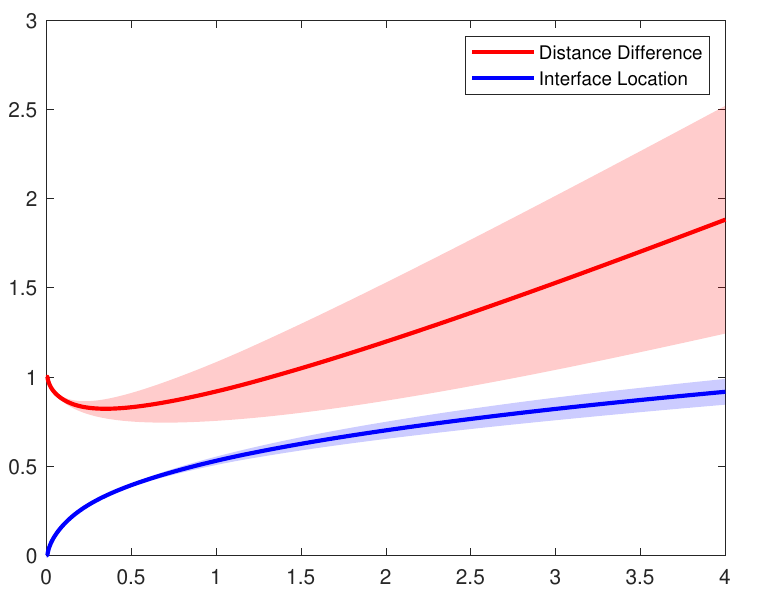}
    \end{minipage}
  }
  \subfigure[$ \widehat{\eta} \sim \mathcal{U}(1.75,2)$,   $\mathbb{E}\left(S^{\star}(4)\right)\approx 0.799$, $\sigma\approx 0.074$]{
    \begin{minipage}[t]{0.45\textwidth}
      \centering
      \includegraphics[width=5.7cm]{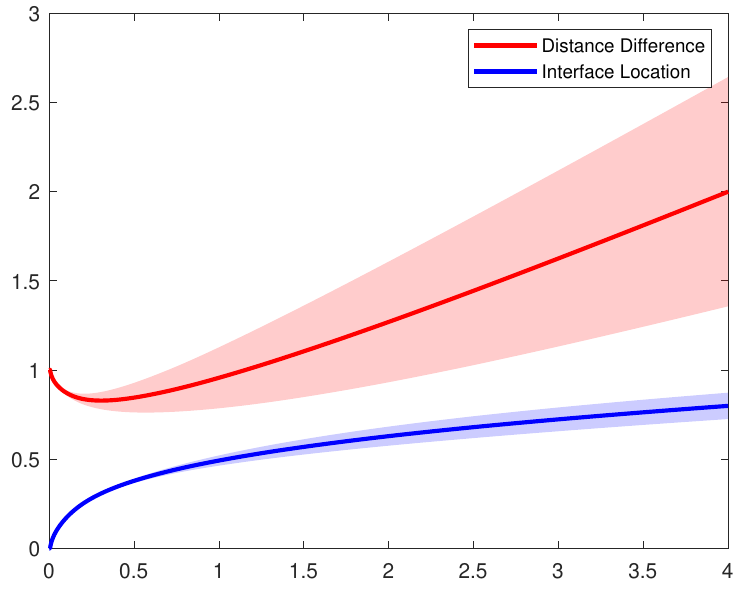}
      \end{minipage}
      
  }
  \caption{
  UQ in in One-dimension: The phase change interface location and distance between inner and outer boundaries are shown for influx uncertainty modeled by $\widehat{\beta} \sim \mathcal{U}(0.2,0.7)$, $T_{0}=-10 ^\circ C$ and different $\widehat{\eta}$ options. Shaded regions indicate plus/minus one standard deviation uncertainty bands on the interface location and boundary separation distance. The mean and standard deviation of the free boundary location $S^{\star}(\tau)$ at $\tau=4$ are also presented.}
       \label{incoming}

  \end{figure}

Two numerical experiments are performed to quantify uncertainty in the phase change interface location. 
In the first one, the airfoil surface temperature $T_{0}$ varies and we choose the IBV $\widehat{\beta} \sim \mathcal{U}(0.2,0.7)$ and incoming heat influx IHI $\widehat{\eta} \sim \mathcal{U}(1,1.25)$. This examines the impact of uncertain airfoil temperature on interface location uncertainty when the IBV and IHI distributions are fixed. 
In the second one, the distribution of the IHI $\widehat{\eta}$ obeys uniform distributions over different ranges while we fix the airfoil temperature $T_{0} = -10$ °C and the IBV distribution $\widehat{\beta} \sim \mathcal{U}(0.2,0.7)$. This investigates the influence of uncertain droplet energy on interface variability when the airfoil temperature and IBV distributions are held fixed.
The mean free interface location and mean difference between the inner interface and injection boundary are shown in Fig. \ref{wingtemp} and Fig. \ref{incoming}, along with $\sigma$ uncertainty bands obtained from generalized polynomial chaos expansions.

The Fig. \ref{wingtemp} depicts the mean interface location and uncertainty bands under different airfoil temperature conditions while keeping the IBV and IHI distributions fixed. It can be observed that the uncertainty in the interface location decreases as the airfoil temperature increases.

The Fig. \ref{incoming} illustrates the results of varying the IHI distribution while fixing the airfoil temperature and IBV distribution. Here, increasing the likely values of the IHI parameter $\widehat{\eta}$ (while keeping interval length fixed) leads to relatively larger uncertainty in the interface location. The uncertainty bands quantify the variability arising from the probabilistic injection parameters in both cases.

The gPC-based uncertainty quantification analysis reveals that random variability in influx parameters significantly impacts phase transition interface behavior. The results highlights the need to account for influx uncertainties when predicting icing regimes. Real-world fluctuations could push systems across rime/glaze thresholds missed by deterministic models. The quantified uncertainty bands provide valuable probabilistic information for developing robust deicing strategies. 

Overall, this analysis advances modeling capabilities to uncover insights difficult to ascertain from deterministic simulations alone. Our analysis highlights the sensitivity of icing regimes to influx fluctuations while providing statistics for robust predictions. Further investigation of the sensitivity mechanisms may elucidate the physical mechanism.
 
\subsection{Uncertainty Quantification in Two-dimension}\label{subsec:UQ2D}
    Now we consider the two-dimensional process with uncertainty. The physical properties and conditions remain the same as Sec. \ref{subsec:Num_2D} except for incoming boundary velocity (IBV)  $\widehat{\beta}$ and incoming heat influx (IHI)  $\widehat{\eta}$.  Fourth-order 2D Legendre polynomials are utilized as the gPC basis functions to construct surrogate models for the interface location. Our model was discretized with $\Delta y=\Delta z=0.005$ and time step $\Delta \tau=0.00001$.

\begin{figure}[ht]
      
  \centering
  \subfigure[$\tau = 0.3$, $y = 0.5$]{
  \begin{minipage}[t]{0.45\textwidth}
  \centering
  \includegraphics[width=6cm]{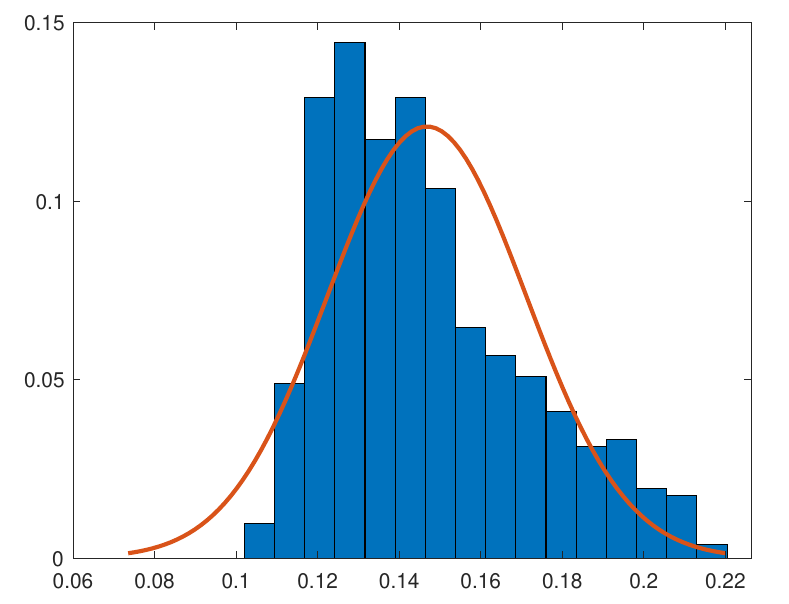}
\label{histogram:a}
  \end{minipage}
  }
  \subfigure[$\tau = 0.7$, $y = 0.5$]{
  \begin{minipage}[t]{0.45\textwidth}
  \centering
  \includegraphics[width=6cm]{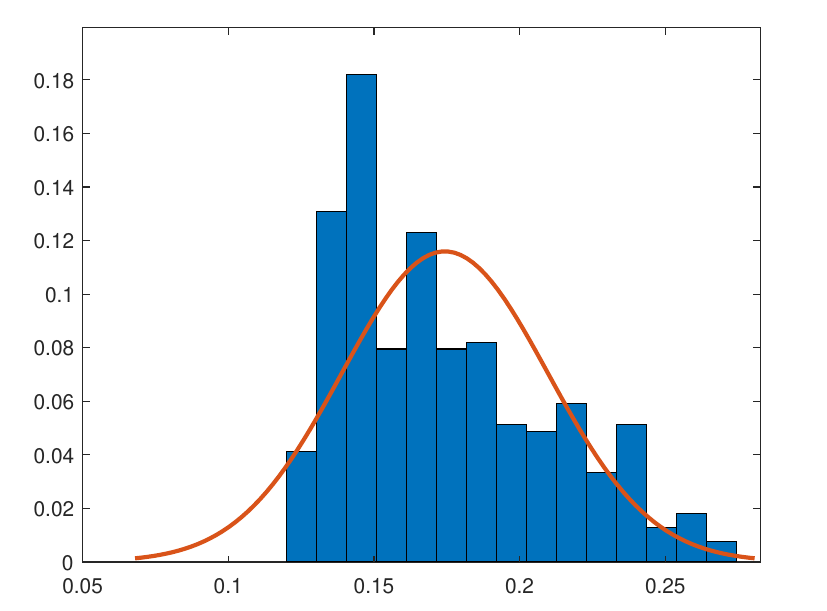}
\label{histogram:b}
  \end{minipage}
  }
  \subfigure[$\tau = 0.3$, $y = 1$]{
  \begin{minipage}[t]{0.45\textwidth}
    \centering
    \includegraphics[width=6cm]{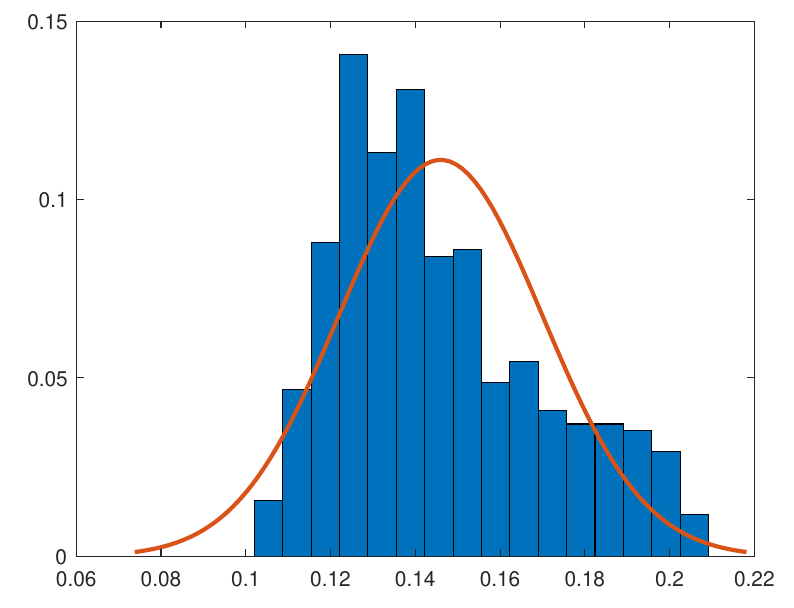}
\label{histogram:c}
    \end{minipage}
  }
  \subfigure[$\tau = 0.7$, $y = 1$]{
    \begin{minipage}[t]{0.45\textwidth}
      \centering
      \includegraphics[width=5.7cm]{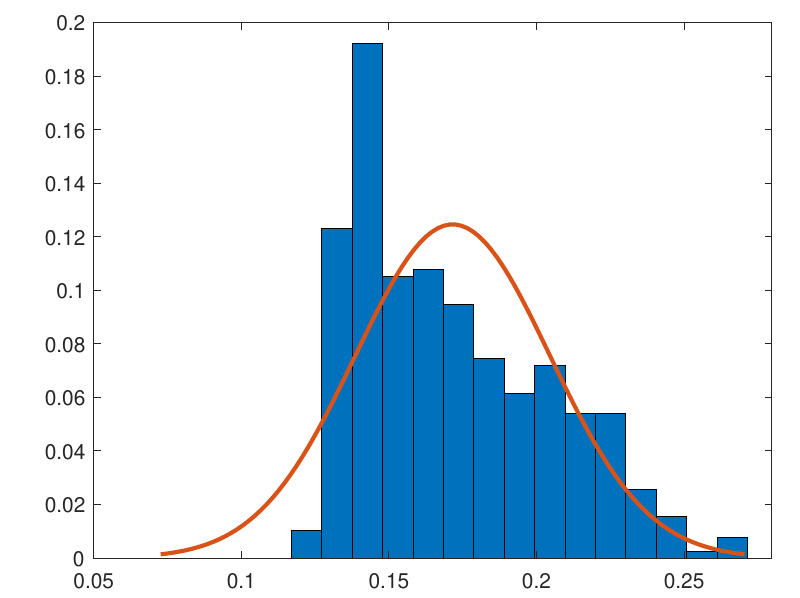}
      \label{histogram:d}
      \end{minipage}
      
  }
  \caption{Histogram distributions of 512 sampled phase interface locations under $\widehat{\beta} \sim \mathcal{U}(0.1,0.3)$ and $\widehat{\eta}= 1+\zeta (1+\cos(3\pi y))$, with $\zeta \sim \mathcal{U}(1,1.2)$ at different spatial positions $y$ and time $\tau$.}
       \label{histogram}

  \end{figure}

     We fixed IBV $\widehat{\beta} \sim \mathcal{U}(0.1,0.3)$ and $\widehat{\eta}= 1+\zeta (1+\cos(3\pi y))$, where $\zeta \sim \mathcal{U}(1,1.2)$. In order to ascertain whether the uncertainty in the 2D free boundary location conforms to a Gaussian distribution across the spatial domain, two sets of numerical experiments are conducted at distinct dimensionless time instances.  The histograms are shown in Fig. \ref{histogram} with 512 sampling points. 
     
  Fig. \ref{histogram:a} and \ref{histogram:c} illustrate that the uncertainty distribution remains relatively stable across different $y$ values when the freezing time is short. Specifically, when $\tau=0.3$, the tight distribution suggests low uncertainty in the output during this early stage. They also highlight the presence of multiple preferred output states influenced by input uncertainties when $\tau=0.3$. On the other hand, Fig. \ref{histogram:b} and \ref{histogram:d} reveal a more pronounced partial peak. As time progresses, the output uncertainty increases, indicating a higher likelihood of encountering extreme situations when $\tau=0.7$.

    \begin{figure}[ht]
      
  \centering
  \subfigure[Skewness of $\tau=0.3$]{
  \begin{minipage}[t]{0.45\textwidth}
  \centering
  \includegraphics[width=6cm]{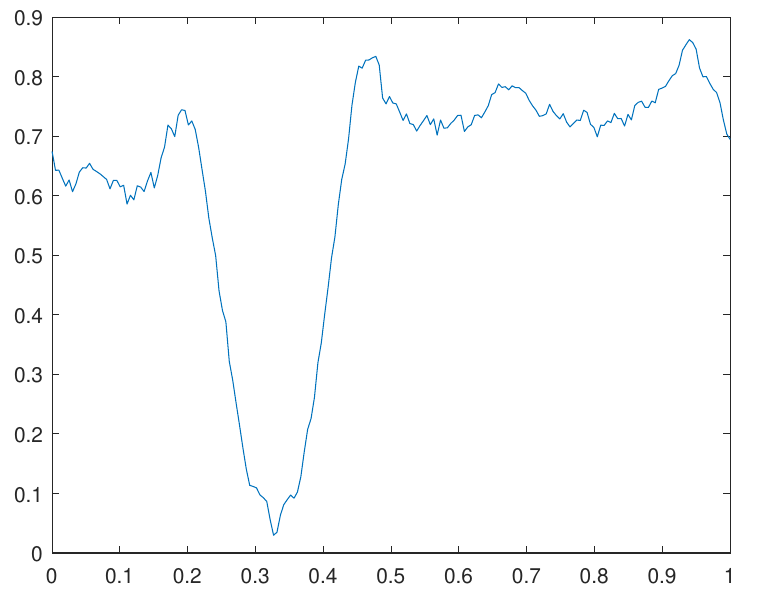}
\label{skew:a}
  \end{minipage}
  }
  \subfigure[Skewness of $\tau=0.7$]{
  \begin{minipage}[t]{0.45\textwidth}
  \centering
  \includegraphics[width=6cm]{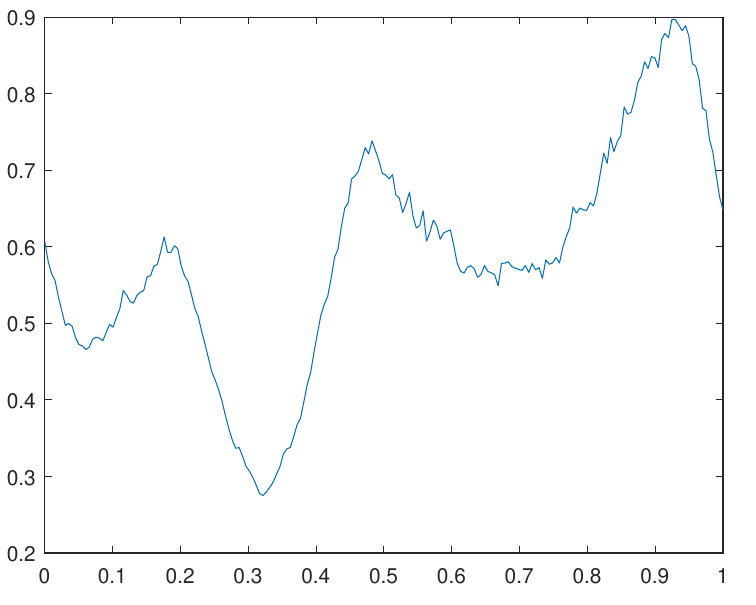}
\label{skew:b}
  \end{minipage}
  }
  \subfigure[Kurtosis of $\tau=0.3$]{
  \begin{minipage}[t]{0.45\textwidth}
    \centering
    \includegraphics[width=6cm]{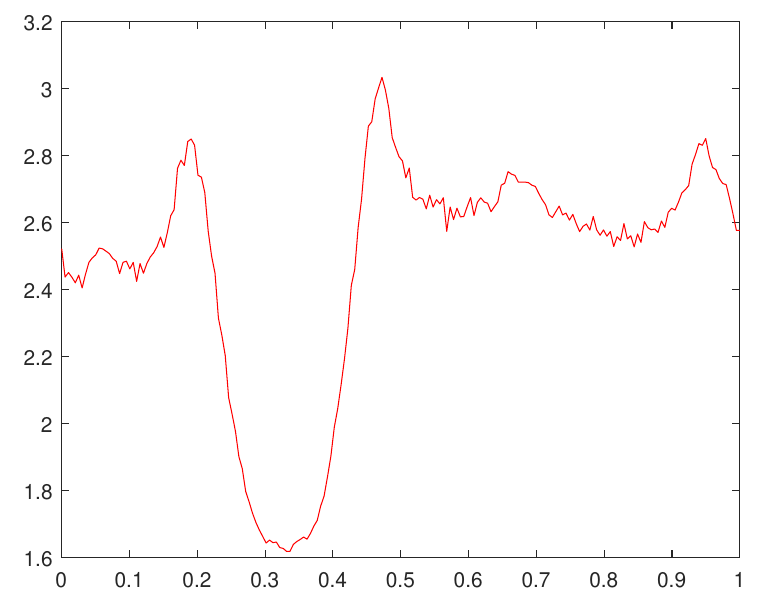}
\label{skew:c}
    \end{minipage}
  }
  \subfigure[Kurtosis of $\tau=0.7$]{
    \begin{minipage}[t]{0.45\textwidth}
      \centering
      \includegraphics[width=5.7cm]{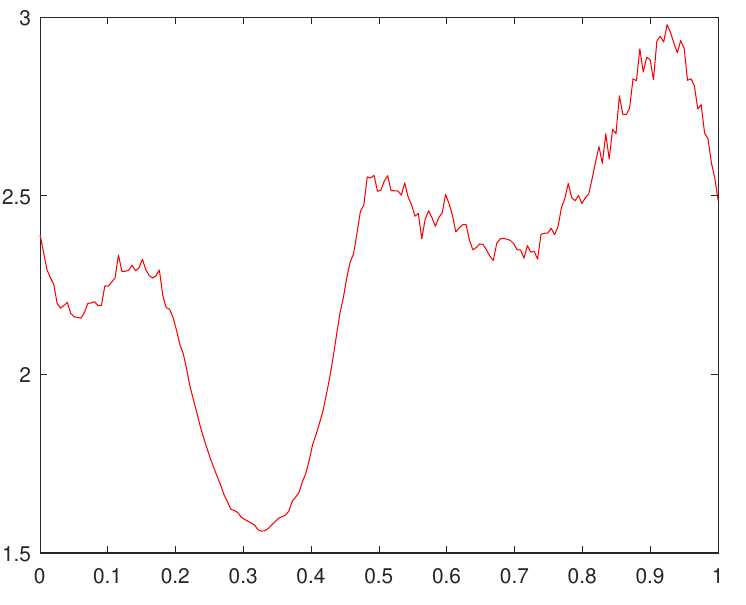}
      \label{skew:d}
      \end{minipage}
      }
  \caption{Skewness and kurtosis of phase interface location distributions under $\widehat{\beta} \sim \mathcal{U}(0.1,0.3)$ and $\widehat{\eta}= 1+\zeta (1+\cos(3\pi y))$, with $\zeta \sim \mathcal{U}(1,1.2)$ across different time $\tau$.}
       \label{skew}

  \end{figure}
    Fig. \ref{skew} illustrates the skewness and kurtosis as functions of $y$ for both $\tau=0.3$ and $\tau=0.7$. When $\tau=0.7$, it reveals a more intricate underlying data distribution, characterized by both leptokurtosis and high skewness. Comparing the data ranges, it is evident that Fig. \ref{skew:b} and \ref{skew:d} exhibit more extreme values than Fig. \ref{skew:a} and \ref{skew:c}, indicating that uncertainty is more likely to result in a complex output scenario when $\tau=0.7$.

 \begin{figure}[ht]
      
  \centering
  \subfigure[$\widehat{\eta}= 1+\zeta+\sin(3\pi y)$, $\zeta \sim \mathcal{U}(1.2,1.5)$]{
  \begin{minipage}[t]{0.45\textwidth}
  \centering
  \includegraphics[width=6cm]{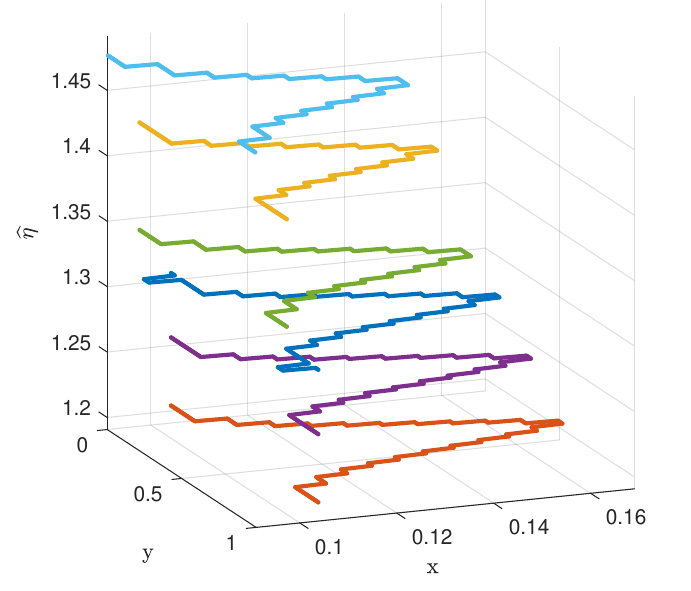}

  \end{minipage}
  }
  \subfigure[$\widehat{\eta}= 2+\sin(3\pi y \zeta )$, $\zeta \sim \mathcal{N}(2,1)$]{
  \begin{minipage}[t]{0.45\textwidth}
  \centering
  \includegraphics[width=6cm]{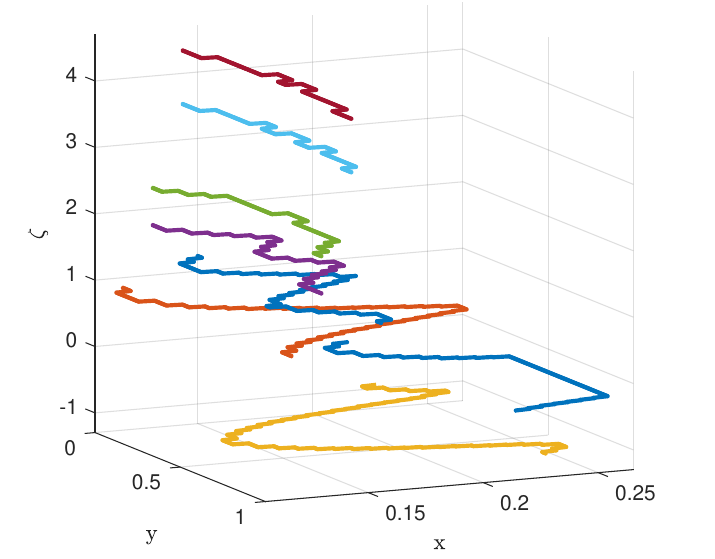}
  \end{minipage}
  }
  \caption{The sampled phase interface location curve in two-dimension under  $\widehat{\beta}=0.2$ and different uncertainty parameters. Each curve which represents the phase change interface location and shape corresponds to a random sample point.}
       \label{final}

  \end{figure}

The results shown in Fig. \ref{final} are also obtained from randomly sampled data points. While 512 sample points are generated, only representative examples are displayed here. Fig. \ref{final} demonstrates that uncertainty in the incoming heat flux has a considerable effect on the morphology and location of the phase transition boundary. This further highlights the importance of accounting for influx uncertainty when predicting freezing states, as deterministic models alone may not provide robust predictions that accurately reflect real-world behavior. Overall, these findings emphasize the need to incorporate probabilistic representations of uncertainty into computational models of phase transition phenomena.

The 2D uncertainty studies unveil non-Gaussian distributions for the free boundary position under uniform variations of parameters. This suggests that complex nonlinear dynamics and internal couplings between multiphysics components may induce non-Gaussian effects even under simple input uncertainties. The plots further demonstrate an escalating complexity in uncertainty propagation through the transient 2D Stefan models. A more in-depth analysis of peak positions, widths, and densities could offer insights into the underlying physics processes. This observation underscores that assumptions of Gaussianity may not always be suitable for uncertainty quantification in multiphysics systems. Characterizing the complete non-Gaussian behavior can provide deeper insights into the model physics.

  \section{Conclusion}\label{sec:Conclu}
In this work, we put forth an integrated modeling and analysis framework for quantifying uncertainties in multiphase systems with coupled injection and phase transition physics. A key innovation is the development of an enthalpy-based model incorporating thermal diffusion across the domain along with energy transfer through an injection boundary. This provides an accurate mathematical representation of the complex dynamics. An implicit enthalpy numerical method is proposed to efficiently handle the interacting moving boundaries and phase transitions. Additionally, customized non-intrusive uncertainty quantification techniques leveraging generalized polynomial chaos and stochastic collocation are adapted to the injection physics context.

Detailed one- and two-dimensional studies demonstrate the integrated framework's ability to effectively capture the competing effects between boundary freezing and water influx under uncertain conditions. Uncertainty quantification of parameters like influx rate and droplet energy content reveals crucial insights into their influence on ice morphology transitions. The higher-dimensional analysis illuminates the subtle dependence between thermodynamic state and injection fluctuations in determining phase transition behavior.

This work considers homogeneous injection boundary growth; an area for future work is to study models with non-homogeneous boundary evolution.  One possibility is leveraging coordinate transformation techniques to handle parametric boundary representations. Challenges include handling topological changes over time. Alternatively, a quasi-3D approach could decompose the domain into subregions with regularized injection shapes. This requires deriving appropriate coupling conditions between subdomains. Dynamic load balancing and high-order geometrical approximation within subregions could help maintain accuracy and computational efficiency as the geometry evolves. In addition, the core modeling and uncertainty quantification techniques established here may generalize to emerging phase-field and variable-order methods for interfacial dynamics. 

   \vskip2mm
\section*{Acknowledgements}
Z. Zhou was supported by the National Key R$\And$D Program
of China, Project Number 2020YFA0712000, and NSFC grant
number 12031013, 12171013. 
The authors are grateful to the anonymous referees for careful reading and
helpful suggestions.

\appendix

\section{Derivation of transformed equations in Sec. \ref{subsec:onecomput}}
\label{Appendix B}
To map the time-dependent domain to a fixed reference domain, we apply the change of variables:
\begin{equation}
    z = \frac{z}{L(t)}, \quad z \in [0,1],
\end{equation}
where $L(t)$ is the length of the growing domain.

With a slight abuse of notation, by writting \(u(x,t)=U(z(x,t),t)\)
where \(u(x,t)\) represents the \(U\) in Eq. (\ref{model:1}),
by chain rule, this gives:
\begin{equation}
  \begin{aligned}
    &\frac{\partial u}{\partial x}=\frac{\partial U}{\partial z}\frac{\partial z}{\partial x}=\frac{\partial U}{\partial z}\frac{1}{L(t)},\\
    &\frac{\partial u}{\partial t}=\frac{\partial U}{\partial t}+\frac{\partial U}{\partial z}\frac{\partial z}{\partial t},\\
    &\frac{\partial z}{\partial t}=-\frac{x}{L^2(t)}L'(t).
\end{aligned}
\end{equation}
Substitute the above relationship into Eq. (\ref{model:1}),
it yields that:
\begin{equation}
  \begin{aligned}
    \frac{\partial U}{\partial t}&=\frac{\partial u}{\partial t}+\frac{\partial U}{\partial z}\frac{x}{L^2(t)}L'(t),\\
    &=\frac{k}{\rho c}\frac{1}{L(t)}\frac{\partial}{\partial z}(\frac{1}{L(t)}\frac{\partial U}{\partial z})+\frac{\partial U}{\partial z}\frac{x}{L^2(t)}L'(t),\\
    &=\frac{k}{\rho c}\frac{1}{L^2(t)}\frac{\partial^2 U}{\partial z^2}+\frac{\partial U}{\partial z}\frac{x}{L^2(t)}L'(t), \quad 0<z<\frac{S(t)}{L(t)}:=S^{\star}(t).
\end{aligned}
\end{equation}
Similarly, we have
\begin{equation}
\begin{aligned}
    \frac{\partial U}{\partial t}=\frac{k}{\rho c}\frac{1}{L^2(t)}\frac{\partial^2 U}{\partial z^2}+\frac{\partial U}{\partial z}\frac{x}{L^2(t)}L'(t), \quad S^{\star}(t)<z<1.
\end{aligned}
\end{equation}
The Stefan condition now reads that:
\begin{equation}
  \begin{aligned}
    \rho L_h\frac{\partial S^{\star}}{\partial t}&=\frac{1}{L(t)}\rho L_h\frac{d S}{d t}+\rho L_h \frac{\partial S^{\star}}{\partial L}\frac{\partial L}{\partial t},\\
    &=\frac{1}{L(t)}(\frac{k}{L(t)}\frac{\partial U}{\partial z}-\frac{k}{L(t)}\frac{\partial U}{\partial z})-\rho L_h\frac{S(t)}{L^2(t)}L'(t),\\
    &=\frac{1}{L(t)}(\frac{k}{L(t)}\frac{\partial U}{\partial z}-\frac{k}{L(t)}\frac{\partial U}{\partial z})-\rho L_h\frac{S^{\star}(t)}{L(t)}L'(t), \quad z=S^{\star}(t).
    \end{aligned}
\end{equation}
Note that stretching and coordinate transformation add one  term to the Stefan condition here.
Substituting these expressions into the original equations Eq. (\ref{eq:category}) yields the transformed system:
\begin{equation}\label{appendixeq:oneConvert}
  \begin{cases}
    \begin{aligned}
        &\rho c\frac{\partial U}{\partial t}=\frac{k}{L^2(t)}\frac{\partial^2 U}{\partial z^2}+\rho c\frac{\partial U}{\partial z}\frac{z}{L(t)}L'(t), \ x \in A, \, B, \ 
 t>0, \\
&U(S^{\star}(t), t)=U(S^{\star}(t), t)=T_m, \quad t>0,\\
& \frac{k}{L(t)} \frac{\partial U}{\partial z}(1,t)+
\frac{d L(t)}{dt} H(1,t)=\eta(t)  \frac{d L(t)}{dt}, \quad t>0 ,\\
&U(z, 0)=T_{\text{initial}}, \quad 0 \leq z \leq 1,\\
&\rho L_h\frac{\partial S^{\star}}{\partial t}=\frac{1}{L(t)}\left[(\frac{k}{L(t)}\frac{\partial U}{\partial z})\right]-\rho L_h\frac{S^{\star}(t)}{L(t)}L'(t),\quad z=S^{\star}(t).  
    \end{aligned}
\end{cases}
\end{equation}

This completes the change of variables mapping the time-dependent domain to the fixed reference domain $z\in [0,1]$.

\section{Proof of the extend enthalpy equation in Section \ref{subsec:onecomput}}
\label{Appendix A.}

Obviously, the Eq. (\ref{eq:gover1}) can not admit a classical pointwise interpretation, since \(H\)
is general not continuous.
Following an usual procedure, we obtain the
weak formulation of Eq. (\ref{eq:gover1}) on multiplying both sides of it by a testing function
\(\Phi \in C_{0}^{\infty}(Q_T)\) and integrating by parts. Here \(d=1\), then
we obtain
\begin{equation}
  \iint_{Q_T}(-H\Phi_t+\frac{k}{L^2(t)}U_z \cdot \Phi_z+\frac{L'(t)}{L(t)}H(z \Phi_z+\Phi))dz dt=0.
\end{equation}
This formulation requires only we give a meaning to the first spatial
derivatives of \(U\). The rigorous definition of the weak formulation
can refer to (cf. \cite{LecStefan,weakformu1,weakformu2}, etc.).
It can be shown this weak formulation admits a solution that also solves the original transformed system Eq. (\ref{eq:oneConvert}). This establishes the enthalpy equation as a proper mathematical representation.

Now we assume \(H\) is a weak solution of Eq. (\ref{eq:gover1}), for any smooth
\(\Phi\) whose support is contained in the interior of 
\(Q_T\).
Considering the regularity of $H$, we can first split the integration into two parts, where in the second part, $H$ is continuous with respect to time $t$.
$$
\begin{aligned}
    \iint_{Q_T}H\Phi_t dzdt&=\int_{0}^{1}(\int_{0}^{T}H\Phi_t dt)dz, \\
    &=\int_{0}^{S^{\star}(T)}(\int_{0}^{T}H\Phi_t dt)dz+\int_{S^{\star}(T)}^{1}(\int_{0}^{T}H\Phi_t dt)dz,
\end{aligned}
$$
where we apply the Fubini theorem. For the first part, we also need to consider its regularity and split it into different parts. For the second part, we can directly adopt the integration by parts formula since $H$ is continuous with respect to  $t$. By applying the integration by parts formula, we have
$$
\begin{aligned}
\int_{0}^{S^{\star}(T)}(\int_{0}^{T}H\Phi_t dt)dz&=
\int_{0}^{S^{\star}(T)}(\int_{0}^{S^{\star -1}(z-)}H\Phi_t dt)dz+\int_{0}^{S^{\star}(T)}(\int_{S^{\star -1}(z-)}^{T}H\Phi_t dt)dz,\\
&=\int_{0}^{S^{\star}(T)} \left(H(S^{\star -1}(z-),z) -H(S^{\star -1}(z+),z)\right) \Phi(S^{\star -1}(z),z) dz\\
&-\int_{0}^{S^{\star}(T)}(\int_{0}^{T}c \rho U_t\Phi dt)dz,
\end{aligned}
$$
since $H_t=c\rho U_t$ both in $A$ and $B$.
By the definition of $H$, the above equation can be simplified as:
$$
\begin{aligned}
\int_{0}^{S^{\star}(T)}(\int_{0}^{T}H\Phi_t dt)dz&=\int_{0}^{S^{\star}(T)} \rho L_h \Phi(S^{\star -1}(z),z)  dz - \int_{0}^{S^{\star}(T)}(\int_{0}^{T}c \rho U_t\Phi dt)dz,\\
&=\int_{0}^{T}\rho L_h \frac{d S^{\star}(t)}{d t} \Phi(t,S^{\star}(t)) dt - \int_{0}^{S^{\star}(T)}(\int_{0}^{T}c \rho U_t\Phi dt)dz,
\end{aligned}
$$
where in the second equation, we introduce a variable substitution $t=S^{\star -1}(z)$.
For the second part, it yields that,
$$
\begin{aligned}
\int_{S^{\star}(T)}^{1}(\int_{0}^{T}H\Phi_t dt)dz=-\int_{S^{\star}(T)}^{1}(\int_{0}^{T}H_t\Phi dt)dz
=-\int_{S^{\star}(T)}^{1}(\int_{0}^{T}c \rho U_t\Phi dt)dz.
\end{aligned}
$$
On adding these two equations, we find
$$
\begin{aligned}
\iint_{Q_T}H\Phi_t dzdt&=\int_{0}^{S^{\star}(T)}(\int_{0}^{T}H\Phi_t dt)dz+\int_{S^{\star}(T)}^{1}(\int_{0}^{T}H\Phi_t dt)dz, \\
&=\int_{0}^{T}\rho L_h \frac{d S^{\star}(t)}{d t} \Phi(t,S^{\star}(t)) dt-\iint_{Q_T}c \rho U_t \Phi dzdt.
\end{aligned}
$$
Then we consider the space part of the differential operator, and it is treated similarly,
$$
\begin{aligned}
 \iint_{Q_T}\frac{k}{L^2(t)}U_z \cdot \Phi_z dzdt&=
 \int_{0}^{T}(\int_{0}^{1} \frac{k}{L^2(t)}U_z \cdot \Phi_z  dz)dt,
\\
&=\int_{0}^{T}(\int_{0}^{S^{\star}(t)} \frac{k}{L^2(t)}U_z \cdot \Phi_z  dz)dt-\int_{0}^{T}(\int_{S^{\star}(t)}^{1} \frac{k}{L^2(t)}U_z \cdot \Phi_z  dz)dt.
\end{aligned}
$$
By applying the integration by parts formula, it reads that
$$
\begin{aligned}
\iint_{Q_T}\frac{k}{L^2(t)}U_z \cdot \Phi_z dzdt&=
\int_{0}^{T}\frac{k}{L^2(t)}\left(U_z^{s} - U_z^{l}\right) \Phi(t,S^{\star}(t)) dt
-\iint_{Q_T} \frac{k}{L^2(t)}U_{zz} \cdot \Phi dz dt,
\end{aligned}
$$
where we denote by $U_z^{s}$ and $U_z^{l}$ the spatial gradient of the restriction of $U$ to solid(liquid).
Again, for the last part,
$$
\begin{aligned}
\iint_{Q_T}\frac{L'(t)}{L(t)}H(z\Phi_z+\Phi) dz dt &=\int_{0}^{T}\frac{L'(t)}{L(t)}\left(\int_{0}^{1}H(z\Phi_z+\Phi )dz\right)dt,\\
=\int_{0}^{T}\frac{L'(t)}{L(t)}\left(\int_{0}^{S^{\star}(t)}H(z\Phi_z+\Phi )dz\right)dt&+\int_{0}^{T}\frac{L'(t)}{L(t)}\left(\int_{S^{\star}(t)}^{1}H(z\Phi_z+\Phi )dz\right)dt.
\end{aligned}
$$
By applying the integration by parts formula, we obatin that
$$
\begin{aligned}
\iint_{Q_T}\frac{L'(t)}{L(t)}H(z\Phi_z+\Phi) dz dt &=\int_{0}^{T}\frac{L'(t)}{L(t)}S^{\star}(t)\left(H(t,S^{\star}(t-))-H(t,S^{\star}(t+))\right)dt \\
&-\iint_{Q_T} \frac{L'(t)}{L(t)} c \rho U_z \Phi dz dt, \\
&=-\int_{0}^{T}\frac{L'(t)}{L(t)}S^{\star}(t) \rho L_h dt-\iint_{Q_T} \frac{L'(t)}{L(t)} c \rho U_z z \Phi dz dt,
\end{aligned}
$$
since $H_z=c\rho U_z$ both in $A$ and $B$.
Thus by combing the above equations
we have,
$$
\begin{aligned}
  &\iint_{Q_T}(-H\Phi_t+\frac{k}{L^2(t)}U_z \cdot \Phi_z+\frac{L'(t)}{L(t)}H(z \Phi_z+\Phi))dz dt\\
  &=\iint_{Q_T}\left(c \rho U_t - \frac{k}{L^2(t)}U_{zz} - c \rho \frac{L'(t)}{L(t)} U_z z \right)\Phi dz dt \\
  &+\int_{0}^{T}\left(-\rho L_h \frac{d S^{\star}(t)}{d t}+ \frac{k}{L^2(t)}\left(U_z^{s} - U_z^{l}\right) -\frac{L'(t)}{L(t)}S^{\star}(t) \rho L_h  \right)\Phi(t,S^{\star}(t)) dt\\
  &=0.
  \end{aligned}
$$
Taking an arbitrary smooth $\Phi$ supported in $A$ and $B$ respectively, then we can derive our original equation,
\begin{equation}
  \rho c\frac{\partial U}{\partial t}=\frac{k}{L^2(t)}\frac{\partial^2 U}{\partial z^2}+\rho c\frac{\partial U}{\partial z}\frac{z}{L(t)}L'(t).
\end{equation}
And the new Stefan condition can be written as
\begin{equation}
\rho L_h\frac{\partial S^{\star}}{\partial t}=\frac{1}{L(t)}\left[(\frac{k}{L(t)}\frac{\partial U}{\partial z})\right]-\rho L_h\frac{S^{\star}(t)}{L(t)}L'(t).
\end{equation}

Then we obtain the solution of Eq. (\ref{eq:oneConvert}), the proof has 
been completed.
\bibliographystyle{siamplain}
\bibliography{UQ_icing}

\end{document}